\documentclass[10pt]{article}
\usepackage{mathrsfs}
\usepackage[tbtags]{amsmath}
\usepackage{epsfig,amstext,amsthm,latexsym}
\usepackage{amssymb,color}
\usepackage{rotating}
\usepackage{multirow}
\usepackage{relsize}
 \allowdisplaybreaks[4]
\usepackage{graphicx, subfigure}
\usepackage[T1]{fontenc}
\usepackage[utf8]{inputenc}
\usepackage{tabularx,ragged2e,booktabs,caption}
\newcolumntype{C}[1]{>{\Centering}m{#1}}

\definecolor{c20}{rgb}{0.,0.7,0.}
\definecolor{c30}{rgb}{0.,0.,1.}
\definecolor{c40}{rgb}{1,0.1,0.7}
\definecolor{c50}{rgb}{1,0,0}

\def\EE#1{\textcolor{c20}{#1}}
\def\EE#1{#1}

\def\ck#1{\textcolor{c20}{#1}}
\def\ck#1{#1}

\def\zz#1{\textcolor{c20}{#1}}
\def\zz#1{#1}

\def\yy#1{\textcolor{c20}{#1}}
\def\yy#1{#1}

\def\xx#1{\textcolor{c20}{#1}}
\def\xx#1{#1}

\def\eee#1{\textcolor{c20}{#1}}
\def\eee#1{#1}

\def\cH#1{\textcolor{c30}{#1}}
\def\cH#1{#1}
\def\aH#1{\textcolor{c30}{#1}}
\def\aH#1{#1}

\def\cw#1{\textcolor{c20}{#1}}
\def\cw#1{#1}

\def\cl#1{\textcolor{c50}{#1}}
\def\cl#1{#1}

\def\cll#1{\textcolor{c50}{#1}}
\def\cll#1{#1}

\def\ccw#1{\textcolor{c40}{#1}}
\def\ccw#1{#1}
\def\peng#1{\textcolor{c40}{#1}}
\def\peng#1{#1}
\def\cL#1{\textcolor{c40}{#1}}
\def\cL#1{#1}
\def\pzx#1{\textcolor{c50}{#1}}
\def\lcx#1{\textcolor{c40}{#1}}
\def\pzx#1{#1}
\def\lcx#1{#1}

\catcode`\@=11 \@addtoreset{equation}{section} \catcode`\@=12

\allowdisplaybreaks[4] 
\newcommand{\nwc}{\newcommand}
\nwc{\COM}[1]{}

\nwc{\vs}[1]{\vskip #1 cm}

\newtheorem{theo}{Theorem}[section]
\newtheorem{sat}[theo]{Proposition}
\newtheorem{de}[theo]{Definition}
\newtheorem{lem}[theo]{Lemma}

\newtheorem{korr}[theo]{Corollary}
\newtheorem{remark}[theo]{Remark}
\newtheorem{exxa}[theo]{Example}

\newcommand{\nelem}[1]{{Lemma \ref{#1}}}

\newcommand{\netheo}[1]{{Theorem \ref{#1}}}
\newcommand{\nekorr}[1]{{Corollary \ref{#1}}}

\def\FRE{\mbox{Fr\'{e}chet }}
\def\d{\mathrm{d}}
\def\I#1{\mathbb{I}\{#1\}}

\def\bar{\overline}

\newcommand{\pk}[1]{{\mathbb{P}} (#1) }

\newcommand{\R}{\!I\!\!R}

\newcommand{\limit}[1]{\lim_{#1 \to   \infty}}

\newcommand{\equaldis}{\stackrel{d}{=}}

\newcommand{\BQN}{\begin{eqnarray}}
\newcommand{\EQN}{\end{eqnarray}}
\newcommand{\BQNY}{\begin{eqnarray*}}
\newcommand{\EQNY}{\end{eqnarray*}}
\newcommand{\BS}{\begin{sat}}
\newcommand{\ES}{\end{sat}}
\newcommand{\BL}{\begin{lem}}
\newcommand{\EL}{\end{lem}}
\newcommand{\BT}{\begin{theo}}
\newcommand{\ET}{\end{theo}}
\newcommand{\BK}{\begin{korr}}
\newcommand{\EK}{\end{korr}}
\newcommand{\BD}{\begin{de}}
\newcommand{\ED}{\end{de}}
\newcommand{\BR}{\begin{remark}}
\newcommand{\ER}{\end{remark}}
\newcommand{\BIT}{\begin{itemize}}
\newcommand{\EIT}{\end{itemize}}
\newcommand{\BDI}{\begin{description}}
\newcommand{\EDI}{\end{description}}

\newcommand{\BEX}{\begin{exxa}}
\newcommand{\EEX}{\end{exxa}}

\newcommand{\QED}{\hfill $\Box$}
\newcommand{\IF}{\infty}

\newcommand{\prooftheo}[1]{ \textsc{Proof of Theorem} \ref{#1} }

\newcommand{\prooflem}[1]{\textsc{Proof of Lemma} \ref{#1}}
\newcommand{\proofkorr}[1]{\textsc{Proof of Corollary} \ref{#1}}

\newcommand{\E}[1]{\mathbb{E}\left\{#1\right\}}

\def\PF{q_t(s)}
\def\PFF{\xx{\Theta_t(s)}}
\begin{document}
%%%%%%%%%%%%%%%%%%%%%%%%%%%%%%%%%%%%%%%%%%%%%%%%%%%%%%%%%
\title{\pzx{Second-order} Tail Asymptotics of Deflated Risks}
\author{\small{Enkelejd Hashorva$^a$\thanks{Corresponding author. Email: enkelejd.hashorva@unil.ch}\quad Chengxiu Ling$^a$\quad Zuoxiang
Peng$^b$}\\
\small{$^a$Faculty of Business and Economics, University of
Lausanne, 1015 Lausanne, Switzerland}\\
\small{$^b$School of Mathematics and Statistics, Southwest
    University, 400715 Chongqing, China}}
\maketitle
\begin{quote}
 {\bf Abstract:}  \pzx{Random deflated risk models have been considered in recent literatures. In this paper, we investigate second-order tail behavior of the deflated risk
  $X=RS$ under the assumptions of second-order regular variation on the survival functions of the risk $R$ and the deflator $S$. Our findings are applied to
approximation of Value at Risk, estimation of small tail probability
under random deflation and tail asymptotics of aggregated deflated
risk}.


\pzx{{\bf JEL classification:} G22}

\pzx{{\bf MSC:} 60G70, 62G32, 62G20, 91B30 }
\end{quote}
%variation; max-domain of attraction; \COM{ ruin
%probability; }  estimation of tail probability.

\section{Introduction}\label{sec1}
 Let $R$ be a non-negative random variable (rv) with distribution function (df) $F$ \pzx{\EE{being} independent of \EE{the} rv $S\in (0,1)$ with df  $G$}.  If $R$ \pzx{models} \ck{the loss amount of a} financial risk, and $S$ \EE{models} a random deflator for a particular time-period, \pzx{then} \EE{the} product  $X=RS$
 represents the deflated value of $R$ at the end of the time-period under consideration. Random deflation is a natural phenomena in most of actuarial applications attributed to the time-value of money. When large values or extremes are of interest, for instance for reinsurance pricing and \ck{risk management purposes},
 it is important to link the \cll{behaviors} of the risk $R$ and the \EE{random} deflator $S$. \ck{Intuitively}, we expect that large values \ck{observed for $R$} are not \ck{significantly} influenced by the \ck{random} deflation.
 However, this is not always the case; a precise analysis driven
 by some extreme value theory models is given in Tang and Tsitsiashvili (2004),  Tang (2006, 2008), %\pzx{Hashorva and Pakes (2010),}
 Hashorva et al.\ (2010), Hashorva (2013), Yang and Hashorva (\cH{2013})\pzx{, Yang and Wang\,(2013), Tang and Yang\,(2012), Zhu and
 Li\,(2012) and the references therein.} \cll{The results} of the
aforementioned papers are obtained mainly under \EE{a}
first\cll{-}order asymptotic condition \pzx{\EE{for the}  survival
function or \EE{the} quantile function in} extreme value theory,
i.e., the df $F$ % \pzx{considered}
\lcx{under consideration} belongs
to the \ck{max-domain of} attraction (MDA) of a univariate extreme
value distribution $\xx{Q_\gamma}, \ck{\gamma\in\R}$,
\ck{abbreviated as} $F\in D(\xx{Q_\gamma})$, \pzx{which means}
\begin{equation}\label{GEV}
F^n(a_nx+b_n)\rightarrow \xx{Q_\gamma}(x):=\exp\left(-(1+\gamma
x)^{-1/\gamma}\right),\quad 1+\gamma x>0,\quad n\to\IF
\end{equation}
 holds for some constants
$a_n>0$ and $b_n\in\R, n\ge1$,  see Resnick\,(1987). The parameter
$\gamma$ is called the extreme value index;  \ck{according} to
$\gamma>0, \gamma=0$ and $\gamma<0$, the df $F$ belongs to the MDA
of the \FRE distribution, the Gumbel distribution and the Weibull
distribution\pzx{, respectively}.

\pzx{In order to derive \EE{some} more informative asymptotic
\EE{results}, second-order regular variation (2RV) conditions are
widely used in extreme value theory. Here we only mention de Haan
and Resnick\,(1996) for the uniform convergence rate of $F^n$ to its
ultimate extreme value distribution $Q_\gamma$ under 2RV, and both
Beirlant et al.\,(2009, 2011), Ling et al.\,(2012) and \lcx{the}
references therein for the asymptotic distributions of
\COM{considered extreme value index estimators.} \lcx{the extreme
value index estimators under consideration.}}

Indeed, \pzx{almost all the common loss distributions
\EE{including}} log-gamma, \EE{absolute} $t$, log-normal, Weibull,
Benktander II, \EE{Beta}, (cf. Table \ref{T1} \peng{in the
Appendix}) \EE{possess 2RV properties}; \EE{actuarial applications
based on those properties are developed} in the recent contributions
Hua and Joe\,(2011), Mao and Hu\,(2012a, \lcx{2012b}) and
Yang\,(2012).

The main contributions of this paper concern \cll{the}
\pzx{second-order expansions of} the tail probability of \pzx{the
deflated risk $X=RS$} \eee{which are then} illustrated by
\ck{several} examples. \eee{Our main findings are} utilized for the
formulations of \cl{three} applications, namely \COM{in the
investigation of asymptotics of ruin probability in some discrete
risk model,}approximation of Value-at-Risk, estimation of small tail
probability \ck{of the deflated risk}, and the derivation of the
tail asymptotics of aggregated risk \ck{under deflation}.

The rest of this paper is organized as follows. Section \ref{sec3}
\ccw{gives} our main results under second-order regular variation
conditions. Section \ref{sec4} \yy{shows} the efficiency of our
second-order asymptotics \yy{through some illustrating examples}.
Section \ref{sec5} is dedicated to \cl{three} applications. The
proofs of all results are \yy{relegated} to Section \ref{sec6}.
\yy{We conclude the paper with a short Appendix.}

\section{Main results}\label{sec3}
\ck{We} start with the definitions and some properties of regular
\eee{variation} \ck{followed by our principal findings.}  A
measurable function $f: \ck{[0, \IF)} \rightarrow \R$ with constant
sign \pzx{near} infinity is said to be of second-order regular
variation with parameters $\alpha\in \R$ and $\rho\le0$, denoted by
$f\in 2RV_{\alpha, \rho}$, if there exists some function $A$ with
constant sign \pzx{near} infinity satisfying
$\lim_{t\rightarrow\infty} A(t)=0$ such that for all $x>0$ (cf.
Bingham et al.\,(1987) and Resnick\,(2007))
\begin{equation}\label{def1}
\lim_{t\rightarrow\infty}\frac{{f(tx)}/{f(t)}-x^\alpha}{A(t)}=x^{\alpha}\int_1^xu^{\rho-1}\,\d
u=:H_{\alpha, \rho}(x).\end{equation}
 Here, $A$ is referred \eee{to}
as the auxiliary function of $f$. \pzx{Noting that \eqref{def1}
implies} $\lim_{t\to\IF}f(tx)/f(t)=x^\alpha$, \pzx{i.e.,} $f$ is
regularly varying at infinity with index $\alpha\in \R$, denoted by
$f\in RV_\alpha$; \eee{$RV_0$ is the class of slowly varying
functions}.
%Clearly, second-order regular variation
%describes the rate of convergence of some function $f\in RV_\alpha,
%\alpha\in \R$.
% Next, we introduce the \cl{definition} of
%second-order $\Pi$-variation.\\
\ck{When $f$} \ck{is} eventually positive, \ck{it is} of
second-order \cl{$\Pi$-}variation with the second-order parameter
$\rho\le0$, denoted by $f\in 2ERV_{0, \rho}$, if there exist some
functions $a$ and $A$ with constant sign \pzx{near} infinity and
$\lim_{t\rightarrow\infty} A(t)=0$ such that for all $x$ positive
\begin{equation}\label{def2}
\lim_{t\rightarrow\infty}\frac{\frac{f(tx)-f(t)}{a(t)}-\log
x}{A(t)}=\psi(x):=\left\{
\begin{array}{cc}
  \frac{x^\rho-1}{\rho}, & \rho<0, \\
  \frac{\log^2x}{2}, & \rho=0
\end{array}
   \right.\end{equation}
\lcx{(cf. Resnick\,(2007))}, where the functions $a$ and $A$ are
referred \yy{to} as the first-order and the second-order auxiliary
functions of $f$, respectively. From Theorem B.3.1 in de Haan and
Ferreira\,(2006), $a\in 2RV_{0,\rho}$ with auxiliary function $A$,
and the second-order auxiliary function $A$ \ck{satisfies}  $|A|\in
RV_{\rho}$. \yy{In fact}, \pzx{\eqref{def2} implies}
$\yy{\limit{t}}(f(tx)-f(t))/a(t)=\log x$ \pzx{for all $x>0$, which
means $f$ is $\Pi$-varying} with auxiliary function $a$, denoted by
$f\in \Pi(a)$.\\
We shall keep the notation of the \cll{introduction
for $R$ and $S\in(0,1)$}, denoting their df's by $F$ and $G$,
respectively, whereas the df of $X=RS$ will be denoted by $H$\pzx{.
Throughout \EE{this paper}, let $\bar F_{0}=1- F_{0}$ denote the survival function
of \EE{a} given distribution $F_{0}$.}

Next, we present our main results in three theorems below.
\netheo{t1} gives a second-order counterpart of \peng{Breiman's
Lemma} (see Breiman\,(1965)) while \netheo{t2} and \netheo{t3}
include refinements of the tail asymptotics of products derived in
Hashorva et al.\,(2010).
\def\eAe{\mathcal{E}}

 \BT\label{t1}
If  $F\in D(Q_{1/\alpha_1})$ satisfies $\bar F\in
2RV_{-\alpha_1,\tau_1}$ with auxiliary function $\tilde{A}$ for some
$\alpha_1>0$ and $\tau_1\le0$, then
\begin{equation}\label{result1}
   \frac{\bar H(x)}{\bar
    F(x)}=\E{S^{\alpha_1}}\left[1+ \eAe(x)\right],
\end{equation}
where
$\eAe(x)=\frac1{\tau_1}\left(\frac{\E\{S^{\alpha_1-\tau_1}\}}{\E{S^{\alpha_1}}}-1\right)\tilde{A}(x)(1+o(1))$
as $x\to \IF$, and thus $\bar H\in 2RV_{-\alpha_1,\tau_1}$ with
auxiliary function
$$A^*(x)=\frac{\E{S^{\alpha_1-\tau_1}}}{\E{S^{\alpha_1}}}\tilde{A}(x).$$
 \ET
{\remark a) The expression for $\tau_1=0$ is understood throughout
this
paper as its limit as $\tau_1\to0$. \\
b) Under the assumptions of \netheo{t1}, \pzx{Breiman's Lemma only}
implies
\begin{align*}
   \frac{\bar H(x)}{\bar
    F(x)}=\E{S^{\alpha_1}}[1+  \eAe^*(x)]
\end{align*}
\pzx{with $\lim_{x\to \IF}  \ck{\eAe^*}(x)=0$, while the} error term
$ \eAe(x)$ in \eqref{result1} not only converges to 0 as $x\to \IF$,
but it shows also the speed of convergence \eee{being} determined by
$\tilde A(x)$.}

\COM{ \label{remark1} Under the conditions of \netheo{t1}, the
second-order tail behavior \ccw{of $X$} is preserved
 under random scaling provided that $\bar F\in 2RV_{-\alpha_1,\tau_1}$. \yy{We have} that $\bar H\in 2RV_{-\alpha_1,\tau_1}$ with auxiliary function
$$A^*(x)=\frac{\E{S^{\alpha_1-\tau_1}}}{\E{S^{\alpha_1}}}\tilde{A}(x).$$
\netheo{t1} provides a possibility to obtain more information on Value-at-Risk,
for more details see Section \ref{VaR_S} below.
}

Next, \eee{we shall consider the \cH{cases} that $F$} belongs to the
MDA of the Gumbel
distribution and the Weibull distribution, respectively.\\
%\xx{For some rv $Y$, with df $Q$ we write sometime $\cH{Y}\sim Q$};
\cll{We write $Y\sim Q$ for some rv \ck{$Y$} with df $Q$,} whereas
$Q^\leftarrow$ denotes the generalized left-continuous inverse of
$Q$ (also for $Q$ which \cH{is} not df). \cll{Since $H$ has the same
upper endpoint \peng{$x_F:=\sup\{y: F(y)<1\}$} as that of the df
$F$,} then all the limit relations below are for $x\to x_F$ unless
otherwise specified. Further, for some $\alpha_2>0$ we set
\begin{equation}\label{Coefficient KK}
L(x)=x^{\alpha_2}\bar G(1-1/x),\quad
K(\alpha_2,\rho)=\left\{\begin{array}{ll}
\frac{(1-\rho)^{-\alpha_2}-1}{\rho}\Gamma(\alpha_2+1), &
\rho<0,\\
\frac{\alpha\cH{_2}\Gamma(\alpha_2+2)}2, &\rho=0, \end{array}\right.
\end{equation}
where $\Gamma(\cdot)$ is the Euler Gamma function, and define
\begin{equation}\label{Coefficient K}
w(x)\eee{=} 1/\E{R-x|R>x},\quad \eta(x)=xw(x).%,\quad x\to \IF.
\end{equation}
\zz{Hereafter the generalized left-continuous inverse of $F$ and $H$
are denoted \cll{by}}
$$U\cl{=}U_R=(1/\bar F)^\leftarrow, \quad U_X=(1/\bar H)^\leftarrow.$$

\BT \label{t2} Let $F$ be strictly increasing and continuous in the
left neighborhood of $x_F$ and \ck{let} $U\in 2ERV_{0,\rho}, \cH{\rho\le0}$
with auxiliary functions $1/w(U)$ and $\tilde{A}$. \ck{If}
\cH{$L\in 2RV_{0, \tau_2}, \tau_2<0$} with auxiliary function $A$,
then
\begin{align}\label{result2}
 \frac{\bar H(x)}{\bar F(x)\bar G\left(1-1/{\eta(x)}\right)}&=\Gamma(\alpha_2+1)+ \eAe(x),
\end{align}
where  $K(\alpha_2,\rho),\eta(x)$ are  defined in \eqref{Coefficient KK},
\cl{\eqref{Coefficient K}} and
$$  \eAe(x)= \Biggl[
     \frac{\Gamma(\alpha_2-\tau_2+1)-\Gamma(\alpha_2+1)}{\tau_2}
   A(\eta(x)) -
   \frac{\alpha_2\Gamma(\alpha_2+2)}{\eta(x)}+K(\alpha_2,\rho)\tilde{A}\left(\frac 1{\bar
    F(x)}\right) \Biggr](1+o(1)).$$
 \ET
\xx{In view of our second result, the error term $ \eAe(x)$ converges
to 0 as $x\to x_F$ with a speed which is determined by $A(\eta(x)),
1/\eta(x)$ and $\tilde A(1/\bar F(x))$. In general, it is not clear
which of these terms is asymptotically relevant for the definition of the error term $\eAe(x)$. For instance in Example
\cH{\ref{Example 2}} below $\tilde A(1/\bar F(x))$ \cH{determines} $ \eAe(x)$. However, Example \cH{\ref{Example
3}} shows the opposite, namely $\tilde A(1/\bar F(x))$ does not appear in
our second-order approximation.}

{\korr \label{Corr1}Under the conditions of \netheo{t2}, \ck{with}  \cH{$\psi$ and $w$ given by \eqref{def2} and} \eqref{Coefficient K},
then for $z\in\R$
\begin{equation}\label{Corr}
   \frac{\bar H(x+z/w(x))}{\exp(-z)\bar H(x)}=1+ \eAe(x), \quad  \eAe(x)= {\left[
   \left(\psi(e^{-z})+\alpha_2\frac{e^{\rho
z}-1}{\rho}\right)
  \widetilde A\left(\frac{1}{\bar
  F(x)}\right)-\frac{\alpha_2z}{\eta(x)}\right]}(1+o(1)),
\end{equation}
where $(e^{\rho z}-1)/{\rho}$ is interpreted as $z$ for $\rho=0$.
\lcx{Thus} $U_X\in 2ERV_{0,0}$ with auxiliary functions $\breve a$
and $\breve A$ given by
\begin{equation}\label{CORR}
\breve{a}(x)=\widetilde a(x)\left(1-\frac{\alpha_2\widetilde
a(x)}{U_X(x)}{+\alpha_2\widetilde A\left(\frac{1}{\bar
F(U_X(x))}\right)}\right),\quad
\breve{A}(x)=-\frac{\alpha_2^2\widetilde{a}^2(x)}{U_X^2(x)}{+
  \widetilde A\left(\frac{1}{\bar F(U_X(x))}\right)},
\end{equation}
where  $\widetilde a=\cll{{1}/{w(U_X)}}$.}

Numerous df's in the MDA of the Gumbel distribution have Weibull
tails (see Embrechts et al.\,(1997) and Table \ref{T2} \cll{in the
Appendix}); specifically such a distribution function $F$ has the
representation \BQN \label{Weib} \bar F(x)=\exp(-V(x)), \EQN where
$V^\leftarrow(x)=x^\theta\ell(x),$ with $\ell$ \peng{denoting} a
positive slowly varying function at infinity.

{\korr \label{Corr2} Under the conditions of \netheo{t2}, if $F$ is
given by \eqref{Weib} \peng{and} $\ell\in 2RV_{0,\rho'},
\cH{\rho'\le0}$ with auxiliary function $b$, then
\begin{align}\label{Corr HF}
\bar H(x)&=\exp(-V(x))\bar
G\left(1-\frac1{V(x)}\right)\Gamma(\alpha_2+1)\theta^{\alpha_2}\left[1+ \eAe(x)
\right],
\end{align}
with
$$  \eAe(x)=\left(\frac{\alpha_2}{\theta}b(V(x))+\frac{\displaystyle\frac{\Gamma(\alpha_2-\tau_2+1)}{\theta^{\tau_2}\Gamma(\alpha_2+1)}-1}{\tau_2}
A(V(x))-\frac{\alpha_2(\alpha_2+1)(\theta+1)}{2V(x)}\right)(1+o(1)),$$
\cll{and thus}
 $$\bar H(x)=\exp(-V^*(x)), \quad (V^*)^\leftarrow(x)=x^\theta\ell^*(x),$$
 where $\ell^*\in 2RV_{0,\max(\rho',-1)}$ with auxiliary function $b^*(x)=b(x)+\theta\alpha_2(\log
 x)/x.$}

\netheo{t1} and \nekorr{Corr1} illustrate that the tail asymptotics
of the product $X=RS$ mainly depends on the heavier factor $R$.
\nekorr{Corr2} shows that for the Weibull tail distributions, the
Weibull tail properties of $X$ are inherited from the factor $R$ in
the presence of random deflation. The result of \nekorr{Corr2} is of
particular interest for the \zz{estimation of tail probabilities},
see Section \ref{ETP}.

\ck{Our last} theorem shows for both $R$ and $S$
belonging to the MDA of the Weibull distribution, the tail of the
product $X=RS$ is heavier than \cll{those} of the factors $R$ and
$S$.

 \BT \label{t3} Let $F$ be strictly increasing and continuous in the left neighborhood of $x_F=1$. Assume
 that for some $\alpha_1>0, \cH{\tau_1\le0}, 1-U\in 2RV_{-1/\alpha_1,\tau_1/\alpha_1}$ with auxiliary function
$\tilde{A}$. \ck{If further}  $L\in 2RV_{0, \tau_2}, \cH{\tau_2\le0}$ with
auxiliary function $A$, then
\begin{align}
 \label{results3}
  \frac{\bar H(x)}{\bar F(x)\bar G(x)}&=\alpha_1B\left(\alpha_1, \alpha_2+1\right)+  \eAe(x),
\end{align}
where
\begin{align*}
   \eAe(x)&=\Biggl[- \frac{\alpha_1^2\alpha_2}{\tau_1}\big(B\left(\alpha_2,\alpha_1-\tau_1+1\right)-B\left(\alpha_2,\alpha_1+1
   \right)\big)
   \tilde{A}\left(\frac1{\bar F(x)}\right)+\frac{\alpha_1}{\tau_2}\big(B\left(\alpha_1,\alpha_2-\tau_2+1\right)
-B\left(\alpha_1,\alpha_2+1\right)\big)A\left(\frac1{1-x}\right)\notag\\
&\quad\quad+\alpha_1\alpha_2B\left(\alpha_1+1,
\alpha_2+1\right)(1-x)\Biggr](1+o(1)),
\end{align*}
with $B(a, b)=\frac{\Gamma(a)\Gamma(b)}{\Gamma(a+b)}, a, b>0$.
 \ET

{\remark \lcx{Recall} that for a df $F$ with a finite endpoint $x_F$
belonging to MDA of the Weibull distribution, then for some
$\alpha_1>0, \tau_1\le0, x_F-U\in 2RV_{-1/\alpha_1,
\tau_1/\alpha_1}$ with auxiliary function $\widetilde A$ is
equivalent to $\bar F(x_F-1/x)\in 2RV_{-\alpha_1, \tau_1}$ with
auxiliary function $\widetilde A^*(x) = -\alpha_1^{2} \widetilde
A\left(\frac1{\bar F(x_F-1/x)}\right)$ and $|\widetilde A^*|\in
RV_{\tau_1}$ (cf. Theorem 2.3.8 in de Haan and Ferreira\,(2006)).
Thus \eqref{results3} holds with
\begin{align*}
   \eAe(x)&=\Biggl[\frac{\alpha_2}{\tau_1}\big(B\left(\alpha_2,\alpha_1-\tau_1+1\right)-B\left(\alpha_2,\alpha_1+1
   \right)\big)
   \tilde{A}^*\left(\frac1{1-x}\right)+\frac{\alpha_1}{\tau_2}\big(B\left(\alpha_1,\alpha_2-\tau_2+1\right)
-B\left(\alpha_1,\alpha_2+1\right)\big)A\left(\frac1{1-x}\right)\notag\\
&\quad\quad+\alpha_1\alpha_2B\left(\alpha_1+1,
\alpha_2+1\right)(1-x)\Biggr](1+o(1)).
\end{align*}
}

 {\remark Under the assumptions of \netheo{t3},
$\bar H(1-1/x)\in2RV_{-\alpha, \tau}$ with
$\alpha=\alpha_1+\alpha_2$ and $ \tau=\max(-1,\tau_1,\tau_2)$.}

\section{Examples}\label{sec4}
In this section, \cl{six} examples are presented to illustrate
estimation errors of the second-order expansions given by Section
\ref{sec3} and the first-order asymptotics by Breiman\,(1965) and
Hashorva et al.\,(2010). We use the \cll{software} R to calculate
the exact value of $\bar H(x)$. 
%Fig.\ \ref{fig1}$\sim$ Fig.\
%\ref{fig3} illustrate \yy{the advantage of our second-order
%tail approximations\cll{.}}
 {\exxa (\FRE case with Pareto
Distribution) Let $R$ be a random variable with a Pareto df $F$
given by
$$\bar F(x)=\left(\frac{\theta}{x+\theta}\right)^\alpha,\quad x>0, \, \alpha, \theta>0$$
denoted in the sequel as  $F\sim Pareto(\alpha, \theta)$. Suppose that $S\sim
beta(a,b)$ \yy{where $beta(a,b)$ stands for the Beta distribution} with positive parameters $a$ and $b$ and probability density function (pdf)
\begin{equation}\label{density of b(a.b)}
g(x)=\frac1{B(a,b)}x^{a-1}(1-x)^{b-1}, \quad 0<x<1,\,  a, b>0.
\end{equation} We have that $\bar
F\in 2RV_{-\alpha,-1}$ with auxiliary function
$\tilde{A}(x)=\alpha\theta/x$ and
$\E{S^\kappa}=B(a+\kappa,b)/B(a,b)$ for all $\kappa>0$. By
\netheo{t1} with \cll{$\alpha_1=\alpha$ and $\tau_1=-1$}
$$\bar H(x)\cH{=}
%\bar F(x)\left(\E{S^\alpha}+\frac{\E{S^{\alpha+1}}-\E{S^\alpha}}{-1}\tilde A(x)\right)=
\bar
F(x)\E{S^\alpha}[1+ \eAe(x)]=\left(\frac{\theta}{x+\theta}\right)^\alpha\frac{B(a+\alpha,b)}{B(a,b)}[1+ \eAe(x)],$$
with
$$ \eAe(x)=\left(1-\frac{\E{S^{\alpha+1}}}{\E{S^{\alpha}}}\right)\tilde{A}(x)(1+o(1))=\frac{\alpha\theta
b}{(\alpha+a+b)x}(1+o(1)).$$ 
%Fig.\ \ref{fig1} compares the
%first-order and the second-order \ck{asymptotic} expansions with the
%exact true value $\bar H(x)$ when \cll{$R\sim Pareto(\alpha,
%\theta), S\sim beta(a,b)$ with $(\alpha, \theta, a, b) = (1,1,1,2)$
%(left) and $(\alpha, \theta, a, b) = (2,1,1,2)$ (right).} \xx{As
%expected, we find that the} second-order tail asymptotics is more
%accurate than the first-order one.
}

{\exxa \label{Example.1} \cl{(\FRE case with Beta distribution of
second kind)} Let $R$ be a random variable with Beta distribution of
second kind with \ck{positive} parameters $a, b$, i.e., $R\stackrel{d}{=} 1/R_0-1,
R_0\sim beta(b,a)$, denoted by $R\sim beta_2(a, b)$
(\cll{here $\stackrel{d}{=}$ stands for equal in distribution}). It
follows from \eqref{density of b(a.b)} that
$$\pk{R_0<x}=\frac{x^b}{bB(b,a)}\left[1-\frac{(a-1)b}{(1+b)}\lcx{x}(1+o(1))\right],\quad x\downarrow0,$$
and thus
\begin{equation}\label{Beta of second order}
\bar
F(x)=\pk{R>x}=\mathbb{P}\left(R_0<\frac1{1+x}\right)=\frac{x^{-b}}{b
B(b, a)}\left[1-\frac{(a+b)b}{(1+b)x}(1+o(1))\right],\quad x\to\IF,
\end{equation}
 i.e., $\bar F\in
2RV_{-b, -1}$ with auxiliary function $\widetilde
A(x)=(a+b)b/((1+b)x)$. \ck{If}  $S\sim beta(c,d)$, then
$\E{S^\kappa}=B(c+\kappa,d)/B(c,d)$ for all $\kappa>0$. \ck{In view of}
\netheo{t1} with $\alpha_1=b$ and $\tau_1=-1$
%, the tail probability
%of $X=RS$, denoted by $\bar H$\cll{, is}
$$\bar H(x)\cH{=}
%\bar F(x)\left(\E{S^\alpha}+\frac{\E{S^{\alpha+1}}-\E{S^\alpha}}{-1}\tilde A(x)\right)=
\bar F(x)\E{S^b}[1+ \eAe(x)]=\frac{x^{-b}}{bB(b,
a)}\left[1-\frac{(a+b)b}{(1+b)x}(1+o(1))\right]\frac{B(c+b,d)}{B(c,
d)}[1+ \eAe(x)],$$ with
$$ \eAe(x)=\left(1-\frac{\E{S^{b+1}}}{\E{S^{b}}}\right)\tilde{A}(x)(1+o(1))=\frac{d}{b+c+d}\frac{(a+b)b}{(1+b)x}(1+o(1)).$$
In particular, for $a=c+d$,
$$\bar H(x)=\frac{x^{-b}}{b B(b, c)}\left[1-\frac{(c+b)b}{(1+b)x}(1+o(1))\right],$$
which is the second-order expansion of survival function of
$beta_2(c, b)$ (cf. \eqref{Beta of second order}), and consistent
with $X\sim beta_2(c, b)$ (see Lemma 5 in Balakrishnan and
Hashorva\,(2011)). 
%Fig.\ \ref{fig1.1} compares the first-order and
%the second-order expansions with the exact true value $\bar H(x)$
%when $R\sim beta_2(a, b), S\sim beta(c, d)$ with $(a, b, c, d)=(3,
%2, 1, 2)$ (left) and $(a, b, c, d)=(2, 2, 1, 2)$ (right). As
%expected, we find that the second-order tail asymptotics is more
%accurate than the first-order one.
}

 {\exxa \label{Example 2}(Gumbel case with
$\rho=0$) Let $R\sim F$ with
\begin{equation}\label{E(1,c)}
\bar F(x)=\exp\left(-\frac{cx}{1-x}\right),\quad 0<x<1, c>0.
\end{equation}
\yy{We write below \eqref{E(1,c)} as $F\sim E(1,c\cL{)}$}. \ck{If follows that}  $F\in D(Q_0)$ with $w(x)= c/{(1-x)^2}$,
and $U\in 2ERV_{0,0}$ with auxiliary functions
$$a(x)\cl{=} 1/w(U(x)), \quad\quad \tilde{A}(x)=-\frac2{c+\log x}.$$ If
$S\sim beta(a,b)$, then we have that the df $G$ of $S$ satisfies
\begin{equation}\label{Beta}
\bar
G\left(1-\frac1x\right)=\frac{x^{-b}}{bB(a,b)}\left(1-\frac{b(a-1)}{(b+1)x}(1+o(1))\right),\quad
x\uparrow\IF,
\end{equation}
i.e., $\bar G(1-1/x)=x^{b}L(x), L\in 2RV_{0,-1}$ with auxiliary
function
$$A(x)=\frac{b(a-1)}{(b+1)\cll{x}}.$$
Consequently
$$\frac1{\eta(x)}=\frac{(1-x)^2}{cx},\quad \widetilde A\left(\frac1{\bar F(x)}\right)=-\frac{2(1-x)}c, \quad A(\eta(x))=\frac{b(a-1)}{(b+1)}\frac{(1-x)^2}{cx}.$$
By \netheo{t2} with $\cll{\alpha_2=b, \tau_2=-1}$ and $\rho=0$,
\begin{align*}
    \bar H(x)&= %\bar F(x)\bar G\left(1-\frac{1}{xw(x)}\right)  \left[\Gamma(b+1)-(\Gamma(b+2)-\Gamma(b+1))
   %A(xw(x))-
   %\frac{b\Gamma(b+2)}{xw(x)}+K(b,0) \tilde{A}\left(\frac 1{\bar
    %F(x)}\right)\right]\\
    %&=
   \bar F(x)\bar G\left(1-\frac{(1-x)^2}{cx}\right) \Gamma(b+1)[1+
   \eAe(x)\pzx{]}
   \end{align*}
   with $$
   \eAe(x)\cH{=}K(b,0)\tilde{A}\left(\frac 1{\bar
    F(x)}\right)(1+o(1))=\frac{b(b+1)}{c}(1-x)(1+o(1)).$$
 }
{\exxa \label{Example 3}(Gumbel case with $\rho<0$) Let $R\sim F$
with
$$\bar F(x)=\frac{1-\exp(-\exp(-x))}{p} ,\quad x>0,\, p=1-e^{-1}.$$
%This function can be treated as the left-truncated Gumbel
%distribution at $x_0=0$.
It follows that  $F\in D(Q_0)$ with constant scaling function
$w(x)=1$ and its tail quantile function is
$$U(x)=\log\frac xp-\frac{p}{2x}(1+o(1)).$$
Furthermore, $U\in 2ERV_{0,-1}$ with auxiliary \ck{functions}
$$a(x)=1,
\quad \tilde{A}(x)=\frac p{2x}.$$ Next, suppose that $S\sim
beta(a,b)$ (cf. \eqref{Beta}\cL{)}. Thus,
$$\frac1{\eta(x)}=\frac1x,\quad \widetilde A\left(\frac1{\bar F(x)}\right)=\frac12e^{-x}, \quad A(\eta(x))=\frac{b(a-1)}{(b+1)x}.$$
By \netheo{t2} with $\cll{\alpha_2=b, \tau_2=-1}$ and $\rho=-1$
\begin{align*}
    \bar H(x)&=
    %\bar F(x)\bar G(1-1/x)\left[\Gamma(b+1)-(\Gamma(b+2)-\Gamma(b+1))
   %A(x)-\frac{b\Gamma(b+2)}x
   %+K(b,-1) \tilde{A}\left(\frac 1{\bar
    %F(x)}\right)\right]\\
    \bar F(x)\bar G\left(1-\frac1x\right) \Gamma(b+1)[1+ \eAe(x)],
\end{align*}
with
$$
 \eAe(x)\cH{=}-\left[\frac{b^2(a-1)}{(b+1)x}+\frac{b(b+1)}{x}\right](1+o(1)).$$
% Fig.\ \ref{fig2} shows that the efficiency of \peng{the} second-order
%asymptotics of $\bar H$ when $R\sim E(1, c)$ with $c=1$ and $S\sim
%beta(1,1/2)$ (left)\peng{;} and $R$ follows the left-truncated
%Gumbel distribution and $S\sim beta(1,1)$ (right).
}

{\exxa \label{Example 3.1}(Gumbel case with Weibull tail) Let $R\sim
\Gamma(\alpha, \lambda)$ with pdf  $$f(x)=
\frac{\lambda^{\alpha}}{\Gamma(\alpha)}x^{\alpha-1}e^{-\lambda
x},\quad x>0, \lambda, \alpha>0.$$
\ck{The  tail} quantile function of $F$ is
$$U(x)\cl{=}\frac{1}{\lambda}(\log x-\log\Gamma(\alpha))\left[1+\frac{(\alpha-1)\log\log x}{\log x-\log\Gamma(\alpha)}(1+o(1))\right].$$
 Thus $F\in D(Q_0)$ with $w(x)=\lambda$ and $U\in 2RV_{0,0}$ with second-order auxiliary function
$$\widetilde A(x)\cl{=}\frac{1-\alpha}{\log^2x}$$ (cf. Table \ref{T2} \peng{in the Appendix}). Next, suppose that
$S\sim beta(a,b)$, where the survival function satisfies
\eqref{Beta}. \ck{Consequently},
$$\frac1{\eta(x)}=\frac1{\lambda x},\quad \widetilde A\left(\frac1{\bar F(x)}\right)=\frac{1-\alpha}{(\lambda x)^2}, \quad A(\eta(x))=\frac{b(a-1)}{(b+1)\lambda x}.$$
By \netheo{t2} with $\cll{\alpha_2=b, \tau_2=-1}$ and $\rho=0$
\begin{align*}
    \bar H(x)&=
    %\bar F(x)\bar G(1-1/x)\left[\Gamma(b+1)-(\Gamma(b+2)-\Gamma(b+1))
   %A(x)-\frac{b\Gamma(b+2)}x
   %+K(b,-1) \tilde{A}\left(\frac 1{\bar
    %F(x)}\right)\right]\\
    \bar F(x)\bar G\left(1-\frac1{\lambda x}\right) \Gamma(b+1)[1+ \eAe(x)],
\end{align*}
with
$$
 \eAe(x)\cH{=}-\frac b{\lambda x}\left[\frac{b(a-1)}{b+1}+(b+1)\right](1+o(1)).$$
Thus
\begin{align}\label{Ex for Weibull}
   \nonumber \bar H(x)&=\frac{(\lambda x)^{\alpha-1}e^{-\lambda
    x}}{\Gamma(\alpha)}\left[1+\frac{\alpha-1}{\lambda x}(1+o(1))\right]\frac{(\lambda
    x)^{-b}\Gamma(b+1)}{bB(a,b)}\left(1-\frac{b(a-1)}{(b+1)\lambda x}(1+o(1))\right)\\
    \nonumber &\quad\cl{\times} \left[1-\frac{b}{\lambda
    x}\left(\frac{b(a-1)}{b+1}+b+1\right)(1+o(1))\right]\\
    &=\frac{(\lambda x)^{\alpha-b-1}e^{-\lambda
    x}}{\Gamma(a)\Gamma(\alpha)/\Gamma(a+b)}\left[1+\frac{\alpha-b(a+b)-1}{\lambda
    x}(1+o(1))\right].
\end{align}
On the other hand, \cl{in view of \nekorr{Corr2}, both $R$ and $X$
are Weibull tail distributions with (cf. Table \ref{T2} \peng{in the
Appendix})
\begin{align}\label{Ex for Weibull2}
\theta=1, \rho'=-1, b(x)=\frac{(1-\alpha)\log x}{x}\quad
\mbox{and}\quad \rho'^*=-1, b^*(x)=b(x)+\frac{\theta \alpha_2\log
x}{x}=\frac{(1-\alpha+b)\log x}{x},
\end{align}
 which is consistent with
\eqref{Ex for Weibull}. In particular, if $\alpha=a+b$, then
\eqref{Ex for Weibull} and \eqref{Ex for Weibull2} are consistent
with the well-known result $X\sim \Gamma(a, \lambda)$ (cf.
Hashorva\,(2013)).}\\
% In Fig.\ \ref{fig2(extra)}, we choose $(\alpha, \lambda, a, b)=(1, 1, 1/2, 1/2)$ (left) and $(\alpha, \lambda, a, b)=(1, 2, 1/2,1/2)$(right).
%We observe that the second-order expansion of the tail probability
%is much closer to the true values.
}

 {\exxa (Weibull case) Let
$R\sim beta(a_1, b_1)$ and $S\sim beta(a_2, b_2)$. By
\COM{Proposition 2.5 in Hua and Joe\,(2011) and} \eqref{Beta},
$1-U\in 2RV_{-1/b_1,-1/b_1}$ with auxiliary function
$$\tilde{A}(x)=-\frac{a_1-1}{b_1(b_1+1)}\left(\frac{x}{b_1B(a_1,b_1)}\right)^{-1/b_1}$$
and $\bar G(1-1/x)=x^{b_2}L(x), L\in 2RV_{0,-1}$ with auxiliary
function
$$A(x)=\frac{b_2(a_2-1)}{(b_2+1)x}.$$
\cll{Hence}
$$\widetilde A\left(\frac1{\bar F(x)}\right)=-\frac{a_1-1}{b_1(b_1+1)}(1-x),\quad A\left(\frac1{1-x}\right)=\frac{b_2(a_2-1)}{b_2+1}(1-x).$$
 By \netheo{t3} with $\alpha_1=b_1,
\alpha_2=b_2, \tau_1=\tau_2=-1$ and
\begin{align*}
\bar H(x)&=%\bar F(x)\bar G(x)\left[b_1B(b_1, b_2+1)
% +b_1^2b_2[B(b_2,b_1+2)-B(b_2,b_1+1)]\tilde{A}\left(\frac1{\bar F(x)}\right)\right.\\
 %&\left.\quad-b_1\left[B\left(b_1,b_2+2\right)
%-B\left(b_1,b_2+1\right)\right]A\left(\frac1{1-x}\right)+b_1b_2B\left(b_1+1,
%b_2+1\right)(1-x)\right]\\
\bar F(x)\bar G(x)\left[b_1B\left(b_1,
b_2+1\right)+ \eAe(x)\right],
\end{align*}
with
$$ \eAe(x)=b_1b_2B(b_1+1,b_2+1)\left(1+\frac{a_1-1}{b_1+1}+\frac{a_2-1}{b_2+1}\right)(1-x)(1+o(1)).$$
\yy{In particular, for}  $a_2+b_2=a_1$,
\begin{align*}
  \bar H(x)&\cll{=}\frac{(1-x)^{b_1+b_2}B(b_1,b_2+1)}{b_2B(a_1,b_1)B(a_2,
  b_2)}\Biggl[1+\Biggl(\frac{b_1+b_2}{b_1+b_2+1}\left(1+\frac{a_1-1}{b_1+1}+\frac{a_2-1}{b_2+1}\right)\\
&\quad-  \left(\frac{b_1(a_1-1)}{b_1+1}+
  \frac{b_2(a_2-1)}{b_2+1}\right)\Biggr)
 (1-x)(1+o(1))\Biggr]\\
  &=\frac{(1-x)^{b_1+b_2}}{(b_1+b_2)B(a_2,b_1+b_2)}\left[1-\frac{(b_1+b_2)(a_2-1)}{b_1+b_2+1}(1-x)(1+o(1))\right],
\end{align*}
 which is the second-order expansion of survival function of
$beta(a_2, b_1+b_2)$ (cf. \eqref{Beta}), and consistent with $X\sim
beta(a_2, b_1+b_2)$ (cf. Hashorva\,(2010)). 
%In Fig.\ \ref{fig3}, we
%simulate the cases with \cll{$(a_1, b_1, a_2, b_2) = (4,2,2,2)$
%(left) and $(a_1, b_1, a_2, b_2) = (\lcx{4},2,2,3)$} (right). \yy{We
%observe that} the second-order expansion of the tail probability is
%much closer to the true values.
}

\COM{

\begin{figure}
\begin{center}
\epsfig{file=Pareto(left).eps, height=150pt, width=200pt,angle=0}
\epsfig{file=Pareto(right).eps, height=150pt, width=200pt,angle=0}
\caption{Tail $\bar H(x)$ when $R\sim Pareto(1,1), S\sim beta(1,2)$
(left) and $R\sim Pareto(2,1), S\sim beta(1,2)$ (right).
 } \label{fig1}
\end{center}
\end{figure}

\begin{figure}
\begin{center}
\epsfig{file=Beta2(left).eps, height=150pt, width=200pt,angle=0}
\epsfig{file=Beta2(right).eps, height=150pt, width=200pt,angle=0}
\caption{Tail $\bar H(x)$ when $R\sim beta_2(3,2), S\sim beta(1,2)$
(left) and $R\sim beta_2(2,2), S\sim beta(1,2)$ (right).
 } \label{fig1.1}
\end{center}
\end{figure}

\begin{figure}
\begin{center}
\epsfig{file=Exponential(c).eps, height=150pt, width=200pt,angle=0}
\epsfig{file=Gumbel(c).eps, height=150pt, width=200pt,angle=0}
\caption{Tail $\bar H(x)$ when $R\sim E(1, c)$ with $c=1$ and $S\sim
beta(1,1/2)$ (left) and $R$ is left-truncated Gumbel distributed and
$S\sim beta(1,1)$ (right).} \label{fig2}
\end{center}
\end{figure}

\begin{figure}
\begin{center}
\epsfig{file=Gamma(left).eps, height=150pt, width=200pt,angle=0}
\epsfig{file=Gamma(right).eps, height=150pt, width=200pt,angle=0}
\caption{Tail $\bar H(x)\sim\Gamma(a, \lambda)$ when $R\sim
\Gamma(\alpha, \lambda)$ and $S\sim beta(a, b)$ for \cl{$(\alpha,
\lambda, a, b)=(1, 1, 1/2, 1/2)$ (left) and $(\alpha, \lambda, a,
b)=(1, 2, 1/2,1/2)$ \cll{(right)}.}} \label{fig2(extra)}
\end{center}
\end{figure}

\begin{figure}
\begin{center}
\epsfig{file=Beta(left).eps, height=150pt, width=200pt,angle=0}
\epsfig{file=Beta(right).eps, height=150pt, width=200pt,angle=0}
\caption{Tail $\bar H(x)$ when $R\sim beta(4,2), S\sim beta(2,2)$
(left) and $R\sim beta(4,2), S\sim beta(2,3)$ (right).
 } \label{fig3}
\end{center}
\end{figure}

}

\newpage

\section{Applications}\label{sec5}
\COM{\subsection{Ruin in the presence of risky investments} Consider
the following insurance risk model. The insurer's net profit
(premium income less claim payment) is denoted by a real-valued
random variable $Z$. And the insurer is in a financial market
consisting of a risk-free bond with a constant interest rate
$\delta>0$ and a risky stock with a stochastic return rate
$\Delta\in(-1,\IF)$. The insurer invests a fraction $\pi\in[0,1)$ of
his current wealth $T_0$ in the stock and keeps the remaining wealth
in the bond. After such an investment, the insurer's wealth is
denoted by $T$ satisfying
$$T=[(1-\pi)(1+\delta)+\pi(1+\Delta)]T_0+Z.$$
The ruin probability $p(x)$ is defined in $p(x)=\pk{T<0|T_0=x}$.
 Assume that $Z$ and $\Delta$ are independent,
denote
\begin{align}\label{Discount factor}&\quad\Upsilon=\frac1{1+\Delta},\quad
R=-Z,\quad S=\frac1{(1-\pi)(1+\delta)+\pi(1+\Delta)}.
\end{align}
Thus $R$ is independent of $S$ and the random variable $\Upsilon$ is
the random discount factor of the risky asset and it takes values in
$(0,\IF)$, the random variable $R$ is the net loss and denote its df
by $F$, while $S$ is \zz{the random discounting factor} \zz{with
distribution function which has upper} endpoint $
x_S=\big[(1-\pi)(1+\delta)\big]^{-1}$.

In view of the first application in Hashorva et al.\,(2010) the
first-order asymptotic behavior of $p(x)$ is
$$p(x)\sim  \bar F(x_1)\pk{\Upsilon> \eee{\eta(x_1)} }\frac{\Gamma(\alpha+1)}{(\pi
x_S)^\alpha}, \quad x_1=x/x_S,$$
where $\eta(x)$ is defined in \eqref{Coefficient K},  provided that $F$ is a
subexponential distribution and $F\in \xx{D(Q_0)}$ with auxiliary
function $w$, and $\pk{\Upsilon>x}\in RV_{-\alpha}$ for some
$\alpha>0$. In our first application we refine the first-order asymptotic expression of $p(x)$.
%Next we give the second-order asymptotics in view of
%\netheo{t2} in the following theorem.
 \BT\label{t_RP} For the risk
model introduced above assume that $F$ is a subexponential
distribution \cl{with $x_F=\IF$} and $U\in 2ERV_{0,\rho},
\cl{\rho\le0}$ with auxiliary functions $1/w(U)$ and $\widetilde A$.
\yy{If for some $\alpha>0, \tau<0$} we have $\pk{\Upsilon>x}\in
2RV_{-\alpha, \tau}$ with auxiliary function $A$ satisfying
\begin{equation}\label{limit a}\lim_{t\to \IF}tA(t)=c\in[-\IF,\IF],
\end{equation}
then
\begin{align*}
    p(x)&\sim \pk{RS>x}= \bar F(x_1)\pk{\Upsilon> \eee{\eta(x_1)}}\frac{\Gamma(\alpha+1)}{(\pi
x_S)^\alpha}[1+ \eAe(x)], \quad x_1=x/x_S,
\end{align*}
with
$$ \eAe(x)=\Biggl[\frac{(\pi x_S)^\tau-1}{\tau}A(\eee{\eta(x_1)})+\frac{\Gamma(\alpha-\tau^*+1)-\Gamma(\alpha+1)}{\tau^*\Gamma(\alpha+1)}
   A^*(\eee{\eta(x_1)})-
        \frac{\alpha(\alpha+1)}{\eee{\eta(x_1)}}+\frac{K(\alpha,\rho)}{\Gamma(\alpha+1)}\tilde{A}\left(\frac 1{\bar
    F(x_1)}\right)\Biggr](1+o(1)),$$
where $K(\alpha,\rho)$ is defined by \cl{\eqref{Coefficient KK}},
and $\tau^*=\max(\tau,-1)$ and auxiliary function
$$A^*(x)=\left\{
                  \begin{array}{ll}
                    -\alpha/x, & c=0, \\
                    (\pi x_S)^\tau A(x), & c=\pm\IF, \\
                    \frac{-\alpha+c/(\pi x_S)}x, & \mbox{otherwise}. \\
                  \end{array}
\right.$$
 \ET
The proof of \netheo{t_RP} is \xx{postponed} to \peng{the}
Appendix.}
\subsection{Asymptotics of \pzx{Value-at-Risk}}\label{VaR_S}
In insurance and \ck{risk management} applications,
\pzx{Value-at-Risk (denoted by $\mathrm{VaR}$)} is an important risk
measure, see e.g., \cl{Denuit et al.\ (2006)}. \ck{We shall analyse
first the} asymptotics of $\mathrm{VaR}_p(X)$ in case that \cll{$R$
has a heavy tail and a Weibull tail, respectively}. Recall that
$\mathrm{VaR}$ at probability level $p$ for some rv $R$ is defined
\peng{by
\begin{equation}\label{VaR and U}
\mathrm{VaR}_p(R)=\inf\{y: F(y)\ge p\}=U(1/(1-p)).
\end{equation}}
 \cll{With the same
notation as before, if $\bar F\in RV_{-\alpha}, \alpha>0$, then by
Breiman's Lemma}
$$\bar H(x)\sim\E{S^\alpha}\bar F(x)\sim \bar F((\E{S^\alpha})^{-1/\alpha}x),\quad x\to\IF$$
implying \peng{the following first-order asymptotics}
\begin{equation}\label{VaR_1}
    \mathrm{VaR}_p(X)\sim (\E{S^\alpha})^{1/\alpha} \mathrm{VaR}_p(R),\quad p\uparrow1.
\end{equation}
\eee{Refining the above, we derive the following second-order asymptotics}
\begin{equation}\label{VaR_2}
    \mathrm{VaR}_p(X)
    = (\E{S^\alpha})^{1/\alpha} \mathrm{ VaR}_p(R)[1+ \eAe(p)], \quad \eAe(p)=\left(\frac{\E{S^{\alpha-\tau}}}{(\E{S^\alpha})^{1-\tau/\alpha}}-1\right)\frac{\tilde{A}(\mathrm{VaR}_p(R))}{\alpha\tau}(1+o(1),\quad p\uparrow1,
\end{equation}
\eee{provided
that $\bar F\in 2RV_{-\alpha,\tau}, \alpha>0, \tau<0$ with auxiliary
function $\tilde{A}$}.\\
  Indeed, there exists some positive constant $c$ such that (cf. Hua
and Joe\,(2011))
$$\bar F(x)=cx^{-\alpha}\left[1+\frac{\tilde{A}(x)}{\tau}(1+o(1))\right]$$
for sufficiently large $x$. Thus, by \netheo{t1}
$$\bar H(x)=cx^{-\alpha}\E{S^\alpha}\left[1+\frac{\E{S^{\alpha-\tau}}}{\E{S^{\alpha}}}\frac{\tilde{A}(x)}{\tau}(1+o(1))\right].$$
Therefore, by Theorem 1.5.12 in Bingham et al.\,(1987) \BQNY
\mathrm{VaR}_p(R)=\left(\frac
c{1-p}\right)^{1/\alpha}\left[1+\frac{\tilde{A}(\mathrm{VaR}_p(R))}{\alpha\tau}(1+o(1))\right],\quad
p\uparrow1 \EQNY and \BQNY \mathrm{VaR}_p(X)&=\left(\frac
{c\E{S^\alpha}}{1-p}\right)^{1/\alpha}\left[1+\frac{\E{S^{\alpha-\tau}}}{\E{S^{\alpha}}}\frac{\tilde{A}(\mathrm{VaR}_p(X))}{\alpha\tau}(1+o(1))\right],\quad
p\uparrow1. \EQNY \eee{Consequently, by} $|\tilde{A}|\in RV_{\tau}$
and \eqref{VaR_1} \eee{we obtain the} second-order asymptotics
\eqref{VaR_2} follows.

\cll{In what follows} we will consider the case that $\bar F$ is in
the MDA of the Gumbel distribution. Since most of such distributions
are Weibull tail distributions (cf. Table \ref{T2} and Tabel
\ref{T1} in the Appendix), we focus on the derivation of the
asymptotics of $\mathrm{VaR}_p(X)$ by $\mathrm{VaR}_p(R)$ (see
\eqref{Varp(XR)} \pzx{below}) under the conditions of
\nekorr{Corr2}. Note that $\bar F$ has a Weibull tail satisfying the
second-order condition (cf. \eqref{Weib})
$$\bar F(x) = \exp(-V(x))\quad \mbox{with}\quad V^\leftarrow(x) =
x^\theta\ell(x)$$ and $\ell\in 2RV_{0,\rho'}, \rho'\le0$ with
auxiliary function $b$. By \eqref{VaR and U}
$$\mathrm{VaR}_p(R) = V^\leftarrow(-\log(1-p)) = (-\log(1-p))^\theta \ell
(-\log(1-p)),\quad p\uparrow1.$$

In view of \nekorr{Corr2} (see \eqref{Corr HF}),
$$\bar H(x) =\exp\left( - V(x) -\alpha_2 \log V(x) + \log L^*(V(x))\right),$$
\pzx{where $L^*$ denotes} a slowly varying function. Recall that
$\log L^*(V(x)) = o(\log V(x))$ (see Bingham et al.\,(1987)), we
have
\begin{align}\label{Varp(XR)}
    \nonumber&\quad \mathrm{VaR}_p(X) = V^\leftarrow\left(-\log(1-p)\left[1-\alpha_2\frac{\log
    (-\log(1-p))}{-\log(1-p)}(1+o(1))\right]\right)\\
    \nonumber &\lcx{=}\left(\log\frac1{1-p}\right)^\theta\left[1-\theta\alpha_2\frac{\log
    \log\frac1{1-p}}{\log\frac1{1-p}}(1+o(1))\right]\ell\left(\log\frac1{1-p}\right) \left[1+\frac{\left(1-\alpha_2\frac{\log
    \log\frac1{1-p}}{\log\frac1{1-p}}\right)^{\rho'}-1}{\rho'}b\left(\log\frac1{1-p}\right)(1+o(1))\right]\\
    & = \mathrm{VaR}_p(R) \left[1-\theta\alpha_2\frac{\log
    \log\frac1{1-p}}{\log\frac1{1-p}}(1+o(1))\right], \quad p\uparrow1.
\end{align}

 \COM{
\subsection{Asymptotics of Conditional Tail Expectation
(CTE)}\label{CTE_S}
 It is well-known that for continuous risks, CTE
is more conservative than VaR and it is defined for a given
probability level $p$ and a rv $R$ by
$$CTE_p(R)=\E{R-u|R>u}, \quad u=VaR_p(R).$$
In this part, we will firstly develop asymptotics of $CTE_p(X)$ by
Hashorva et al.\,(2010). For some $x\in(0,1)\cup(1,\IF)$
\begin{align}\label{CTE_1}\frac{CTE_p(X)}{VaR_p(X)}-1&\sim\frac{{VaR_{p'}(X)}/{VaR_p(X)}-1}{\log
x},\quad p\uparrow1, p'=1-(1-p)/x\end{align} provided that $F\in
D(Q_0)$. Secondly, we will investigate that for $\bar F\in
2RV_{-\alpha,\tau}, \alpha>1,\tau<0$ with auxiliary function
$\tilde{A}$, second-order asymptotics of CTE of $X$ is
\begin{equation}\label{CTE_2}
     \frac{CTE_p(X)}{(\E{S^\alpha})^{1/\alpha}
     CTE_p(R)}\sim \left(1+\aH{\frac{\alpha-1}{\alpha-1-\tau}\left(\frac{\E{S^{\alpha-\tau}}}{(\E{S^\alpha})^{1-\tau/\alpha}}-1\right)}\frac{\tilde{A}(VaR_p(R))}{\alpha\tau}\right),\quad p\uparrow1.
\end{equation}
Now we prove \eqref{CTE_1}. Note that the auxiliary function $w$ for
$F\in D(Q_0)$  is self-neglecting, i.e., $w(t+x/w(t))\sim w(t)$,
implying that random contraction $X$ still belongs to MDA of Gumbel distribution. Thus
\begin{equation}\label{3.1}
\E{R-x|R>x}\sim \frac1{w(x)}\sim \E{X-x|X>x},\quad
\frac{U_X(tx)/U_X(t)-1}{\log x}\sim\frac1{U_X(t)w(U_X(t))}.
\end{equation}
Substitute $t=1/(1-p),
  tx=1/(1-p')$ into \eqref{3.1}, the result in \eqref{CTE_1} follows.
\\
For the claim in \eqref{CTE_2}, by Proposition 5 in Hua and
Joe\,(2011)
$$\frac{CTE_p(R)}{VaR_p(R)}\sim\frac{\alpha}{\alpha-1}\left(1+\frac{\widetilde A(VaR_p(R))}{\alpha(\alpha-1-\tau)}\right), \quad
\frac{CTE_p(X)}{VaR_p(X)}\sim\frac{\alpha}{\alpha-1}\left(1+\frac{
A^*(VaR_p(X))}{\alpha(\alpha-1-\tau)}\right)$$ with $A^*$ defined in
Remark \ref{remark1}. Therefore, \eqref{CTE_2} follows from Remark
\ref{remark1}, \eqref{VaR_2} and straightforward calculations.
}
\subsection{Estimations of tail probability}\label{ETP}
\ck{In many insurance applications it is important to estimate the
tail probability of the extreme risks.} In what follows, we
investigate this problem \ck{under the random scaling framework}.
Let $\{(R_i,S_i), i=1,\cdots, n\}$ \ck{be a random sample} from
$(R,S)$, our goal is to estimate $p=\pk{X>x}=\pk{RS>x}$ with
sufficiently large $x$. One \cll{possible estimation} is via the
empirical df if $x$ is in the region of the sample $X_i,i\le n$ with
$X_i=R_iS_i, i=1,\cdots,n$. In general, we consider how to estimate
$p_n:=\pk{X>x_n}$ \pzx{as} $x_n\to\IF$. Let $R_{n-k+1,n}, S_{n-k+1,
n}$ and $X_{n-k+1,n}, k=1,\cdots, \cll{n}$ be the associated
increasing order statistics and $R\sim F$ and $S\in (0,1)$ are
independent.

First we consider the case that $\bar F\in2RV_{-\alpha,\tau}$ with
$\alpha>0,\tau<0$ and the second-order auxiliary function
$\tilde{A}$, thus \lcx{by Hua and Joe\,(2011), there exists a
positive constant $c$ such that}
$$\bar F(x)\cl{=} c
x^{-\alpha}(1+\widetilde A(x)/\tau(1+o(1)))=:c
x^{-\alpha}(1+\alpha\delta(x)),$$ i.e., $F\in
\mathcal{F}_{1/\alpha,\tau}$ with $\delta(x)=\widetilde
A(x)/(\alpha\tau)$ (cf. Beirlant et al.\,(2009)). By \netheo{t1}
\begin{align}\label{Estimation1}
\bar H(x)\cl{=} %\bar
%F(x)\left(\E{S^{\alpha}}+\mathbb{E}\left\{S^{\alpha}\frac{S^{-\tau}-1}{\tau}\right\}\tilde{A}(x)\right)=
\bar
F(x)\left(\E{S^{\alpha}}+\mathbb{E}\left\{S^{\alpha}(S^{-\tau}-1)\right\}\alpha\delta(x)(1+o(1))\right).
\end{align}
In order to estimate $\bar H(x)$ with $x=x_{n}$, we use the
estimators of $1/\alpha, \delta, \tau$ and $\bar F$ proposed by
Beirlant et al.\,(2009). Let $y_{k,n}=x_n/R_{n-k,n},
\widehat\tau_{k,n}=\widehat\rho_n/H_{k,n}$ with $\widehat\rho_n$
some weakly consistent estimator of $\rho=\tau/\alpha$ based on
samples from the parent $R$, denote
$$H_{k,n}=\frac1k\sum_{i=1}^k\log R_{n-i+1,n}-\log
R_{n-k,n},\quad E_{k,n}(s)=\frac1k\sum_{i=1}^k
(R_{n-i+1,n}/R_{n-k,n})^s,s\le0$$ and
\begin{equation*}
\widehat\alpha_{k,n}=\left(H_{k,n}-\widehat\delta_{k,n}\frac{\widehat\rho_n}{1-\widehat\rho_n}\right)^{-1},\quad
    \widehat\delta_{k,n}=H_{k,n}(1-2\widehat\rho_n)(1-\widehat\rho_n)^3\widehat\rho_n^{-4}\left(E_{k,n}(\widehat\rho_n/H_{k,n})-\frac{1}{1-\widehat\rho_n}\right).
\end{equation*}
Thus, by \eqref{Estimation1}, \cll{the tail probability} $p_n$ can
be estimated as (denoted by $\widehat p_{k,n}(R,S)$)
\begin{equation}\label{Based on R, S}
    \widehat p_{k,n}(R,S)=\widehat{\bar
    F}(x_n)\left(\widehat{\E{S^{\alpha}}}+(\widehat{\E{S^{\alpha-\tau}}}-\widehat{\E{S^{\alpha}}})\frac{\widehat{\delta}_{k,n}}{H_{k,n}}\right),
\end{equation}
with
\begin{align}\label{estimation1}
    \widehat{\bar F}(x_n)&=\frac
    kn\left(y_{k,n}\left(1+\widehat\delta_{k,n}(1-y_{k,n}^{\widehat\tau_{k,n}})\right)\right)^{-\widehat\alpha_{k,n}},\quad \widehat{\E{S^{\alpha}}}=\frac1n\sum_{i=1}^n
S_i^{\widehat\alpha_{k,n}},\quad
\widehat{\E{S^{\alpha-\tau}}}=\frac1n\sum_{i=1}^n
S_i^{\widehat\alpha_{k,n}-\widehat\tau_{k,n}}.
 \end{align}
On the other hand, by \netheo{t1}, $X$ has the same second-order
tail behavior as $R$. Consequently, $p_n$ can be directly estimated
by using samples from $X$. We denote \ck{that estimator} (cf. \eqref{estimation1}) by $\widehat p_{k,n}(X)$, given as
\begin{equation}\label{Based on X}
    \widehat p_{k,n}(X)=\frac
    kn\left(y^*_{k,n}\left(1+\widehat\delta^*_{k,n}(1-(y^*_{k,n})^{\widehat\tau^*_{k,n}})\right)\right)^{-\widehat\alpha^*_{k,n}},
    \end{equation}
    with $y_{k,n}^*=x_n/X_{n-k,n}$ and $\widehat\delta^*_{k,n},
    \widehat\tau^*_{k,n}, \widehat\alpha^*_{k,n}$ are $\widehat\delta_{k,n},
    \widehat\tau_{k,n}, \widehat\alpha_{k,n}$ with the order statistics replaced by $\{X_{n-k+1,n}, k=1,\cdots, n-1\}$.
    \\
    \ck{Relying on} \eqref{Based on R, S} and \eqref{Based on X},
    we \ck{shall} perform some simulations to compare $\widehat\alpha_{k,n}, \widehat p_{k,n}(R, S)$
    and $\widehat\alpha^*_{k,n}, \widehat p_{k,n}(X)$. Since $\tau=-1$ holds in most applications, we take $\widehat\tau_{k,n}=-1$
 and $\widehat\rho_n=-H_{k,n}$. 
 %For each simulation in Fig.\ \ref{fig4}$\sim$ Fig.\ \ref{fig6}, we take one sample from Table \ref{T1} \peng{in the Appendix} with \cll{sample} size % $n=1000$
%  and estimate $1/\alpha$ and $p=\pk{X>3}$.
%
% In Fig.\ \ref{fig4}, let $R\sim Pareto(2,1)$ and $ S\sim beta(1,2)$, thus $\bar F\in 2RV_{-\alpha,\tau}$ with $\alpha=2, \tau=-1$. We estimate %$1/\alpha=0.5$ and $p=0.01298$. Clearly, the
% estimators\cll{,} based on the original data information on $R$ \cll{and $S$,} are asymptotically unbiased and have much wider stable regions
% than those \pzx{directly} based on $X=RS$.

%In Fig.\ \ref{fig5}, let $R$ be \COM{truncated t} \cL{absolute $t$}
%distributed with freedom degree $v=4$ and $S\sim beta(1,2)$, thus
%$\bar F\in 2RV_{-\alpha,\tau}$ with $\alpha=4, \tau=-2$. The
%estimation of $p=0.0035$ based on our theorem is much closer to the
%true values, even in our simulations we take
%$\widehat\tau_{k,n}=-1$.

\COM{
In Fig.\ \ref{fig6}, let $R$ be Log-gamma distributed with shape and
scale parameters $\alpha=4, \beta=2$ and $S\sim beta(2,4)$, thus
$\bar F\in 2RV_{-\alpha,\tau}$ with $\alpha=4, \tau=0$. This does
not satisfy the assumption $\tau<0$ for the  \eee{estimators} in
\eqref{Based on R, S} and \eqref{Based on X}. However, the
estimation of $p=0.0037$ based on our theorem performs reasonably
well.
}
Next, we investigate the case \cl{of} $F\sim D(Q_0)$. For
convenience, we consider \ck{only} the estimation comparisons for $F$
being Weibull tail distributions. By \nekorr{Corr2}, \cll{both} $R$
and $X$ are Weibull tail distributions with the same Weibull tail
coefficient $\theta$ and further the second-order parameter
$\rho'^*$ is greater than $-1$, we consider the bias-reduced Weibull
tail coefficient estimators \peng{$\widehat\theta$ due to Diebolt et
al.\,(2008):}
\begin{align}
    \label{e_thetaR}&\quad\widehat\theta=\widehat\theta(k, R)=\bar Z_k-\widehat b(\log(n/k))
    \bar x_k,
\end{align}
with
$$\widehat b(\log(n/k))=\frac{\sum_{i=1}^k(x_i-\bar x_k)Z_i}{\sum_{i=1}^k(x_i-\bar
    x_k)^2}$$
and
$$x_j=\frac{\log(n/k)}{\log(n/j)}, \quad Z_j=j\log(n/j)(\log
R_{n-j+1,n}-\log
 R_{n-j,n}),\quad \bar
x_k=\frac{\sum_{j=1}^k x_j}{k},  \quad \bar Z_k=\frac{\sum_{j=1}^k
Z_j}{k}.$$ Based on the bias-reduced tail quantile estimators
provided by Diebolt et al.\,(2008), \ck{given} by
$$\widehat x_{p_n}=R_{n-k,n}\left(\frac{\log(1/{p_n})}{\log(n/k)}\right)^{\widehat\theta}\exp\left(\widehat b(\log(n/k))\frac{(\log(1/{p_n})/\log(n/k))^{\widehat\rho'}-1}{\widehat\rho'}\right)$$
with $p_n$ known, we \cll{can} solve the dual problem and estimate
the tail probability $\bar F(x)$ for given $x$ as follows
\begin{equation}\label{e_pR}
    \widehat{\bar F}(x)=\exp\left(-\log(n/k)\left(\frac {x}{R_{n-k,n}}\right)^{1/\widehat\theta}\exp\left(-\widehat b(\log(n/k))\frac{(x/R_{n-k,n})^{\widehat\rho'/\widehat\theta}-1}{\widehat\theta\widehat\rho'}\right)\right),
\end{equation}
with $\widehat\rho'$ a consistent estimator of $\rho'$. Note that
$\bar F(x)=\exp(-V(x))$ and $S\sim G$ with $\bar G(1-1/x)\in
2RV_{-\alpha_2,\tau_2}$ is equivalent that $S^*=1/(1-S)\sim G^*$
with $\bar G^*\in 2RV_{-\alpha_2,\tau_2}$. Hence by \eqref{e_pR} and
Beirlant et al.\,(2009), we have
\begin{equation}\label{e_V}
\widehat V(x)=-\log \widehat{\bar F}(x),\quad \widehat
b(V(x))=\widehat b(\log(n/k))\left(\frac{\widehat
V(x)}{\log(n/k)}\right)^{\widehat\rho'}
\end{equation}
and
\begin{align}\label{e_G}
    \widehat{\bar
    G}(1-1/V(x))=\frac{k}{n}\left(y_{k,n}(1+\widehat\delta_{k, n}(1-y_{k,n}^{\widehat{\tau}_2(k)}\right)^{\widehat{\alpha}_2(k)},
    \quad \widehat A(V(x))=\widehat{\alpha}_2(k)\widehat{\tau}_2(k)\widehat\delta_{k, n}\left(y_{k,n}\right)^{\widehat{\tau}_2(k)},
\end{align}
where $y_{k,n}=\widehat{V}(x)/S^*_{n-k,n}$ and
$\widehat\delta_{k,n}, \widehat{\tau}_2(k), \widehat{\alpha}_2(k)$
are estimated with \cll{the} order statistics replaced by
$S^*_{n-k,n}:=1/(1-S_{n-k,n})$ in \eqref{estimation1}. Therefore,
combining \eqref{e_thetaR}, \eqref{e_pR}, \eqref{e_V} and
\eqref{e_G}, the \yy{estimator} of $p=\bar H(x)$, \lcx{denoted by
$p_k(R, S)$, is given as} (cf. \nekorr{Corr2})
\begin{align}\label{e_pRS}
   \nonumber \widehat p_k(R, S)&=\widehat{\bar F}(x)\widehat{\bar
    G}(1-1/V(x))\Gamma(\widehat{\alpha}_2(k)+1)(\widehat\theta)^{\widehat{\alpha}_2(k)}\\
    &\quad\times \left[1+\frac{\widehat\alpha_2(k)}{\widehat\theta}\widehat
b(V(x))+\frac{\displaystyle\frac{\Gamma(\widehat\alpha_2(k)-\widehat{\tau}_2(k)+1)}{(\widehat\theta)^{\widehat{\tau}_2(k)}\Gamma(\widehat{\alpha}_2(k)+1)}-1}{\widehat{\tau}_2(k)}
\widehat
A(V(x))-\frac{\widehat{\alpha}_2(k)(\widehat{\alpha}_2(k)+1)(\widehat\theta+1)}{2\widehat
V(x)}\right].
\end{align}
On the other hand, by \nekorr{Corr2}, we can estimate $p=\bar H(x)$
directly based on samples from $X$ as
\begin{equation}\label{e_pX}
   \widehat p_k(X)=\exp\left(-\log(n/k)\left(\frac x{X_{n-k,n}}\right)^{1/\widehat\theta^*}\exp\left(-\widehat
   b^*(\log(n/k))\frac{(x/X_{n-k,n})^{\widehat\rho'^*/\widehat\theta^*}-1}{\widehat\theta^*\widehat\rho'^*}\right)\right),
\end{equation}
where $\widehat\rho'^*$ is a consistent estimator of $\rho'^*$ and
$\widehat\theta^*, \widehat b^*$ are computed by \eqref{e_thetaR}
with samples replaced by $X_i=R_iS_i, i=1,2,\cdots,n$.

Now, we perform the simulations of the estimators of $\theta$ and
$p=\pk{X>x}$ given by \eqref{e_pRS} and \eqref{e_pX} with one sample
of size $n=5000$ from Table \ref{T2} and Table \ref{T1} \peng{in the
Appendix}. In the simulation we take $\widehat\tau_{2}(k)=-1,
\widehat\rho'=\widehat\rho'^*=-1$ and plot sample paths of
$\widehat\theta$ and $\log(\widehat p_k/ p), k=100,\cdots,4500$,
\cH{with $\widehat p_k=\widehat p_k(R, S), \widehat p_k(X)$,
respectively (cf. \eqref{e_pRS} and \eqref{e_pX}).}

\COM{ In Fig.\ \ref{fig7}, let $R\sim |N(0,1)|$ and $S\sim beta(1,3)$, thus $\theta=1/2, \rho'=-1$ and $\alpha_2=3, \tau_2=-\IF$. The true
 value of $\pk{X>3}$ is $5.281\times 10^{-6}$. Clearly, the
 \eee{estimators} based on the original data information on $R, S$ are asymptotically unbiased and have much wider stable regions than
 those \yy{based} on
 $X=RS$.

In Fig.\ \ref{fig8}, let $R\sim W(2,1)$ and $S\sim beta(2,3)$, thus
$R$ is sub-exponential distributed with $\theta=1/2, \rho'=-\IF$ and
$\alpha_2=3, \tau_2=-1$. The estimation of $\theta$ and
$\pk{X>3}=2.1186\times 10^{-7}$ based on our theorem is much closer
to the true values, even in our simulations we take
$\widehat\rho'=-1$.

In Fig.\ \ref{fig9}, let $R\sim \Gamma(1/2, 1)$ and $S\sim beta(1,3)$,
thus $R$ is exponential-like distributed with $\theta=1, \rho'=-1$
and $\alpha_2=3, \tau_2=-\IF$. The  \eee{estimators} of $\theta$ and
$\pk{X>10}=1.7189\times 10^{-8}$ based on \eqref{e_pRS} performs
reasonably better than those based on \eqref{e_pX}.

In Fig.\ \ref{fig10}, let \lcx{$R\sim PW(9/10,9/20)$ with $C=1,
D=2$} and $S\sim beta(1,3)$, thus $R$ is super-exponential
distributed with $\theta=10/9, \rho'=-1/2$ and $\alpha_2=3,
\tau_2=-\IF$. The
 \eee{estimators} of $\theta$ and $\pk{X>10}=3.5539\times 10^{-9}$ based
on \eqref{e_pRS} performs reasonably better than those based on
\eqref{e_pX}, even in our simulations we take $\widehat\rho'=-1$ and
 $\widehat\tau_2=-1$ for the true values $\rho'^*=-1/2, \tau_2=-\IF$.
}
\COM{
\begin{figure}
\begin{center}
\subfigure[$R\sim Pareto(2,1),S\sim beta(1,2).$
]{%
            \label{fig4}
 \epsfig{file=Pareto(tp).eps, height=200pt, width=160pt,angle=0}
                   }%
                   \subfigure[$R\sim |t|(4), S\sim beta(1,2).$]{%
            \label{fig5}
 \epsfig{file=Absolute(t).eps, height=200pt, width=160pt,angle=0}
                   }%
                   \subfigure[$R\sim Log\cll{\text{-}}gamma(4,2), S\sim beta(2,4).$]{%
            \label{fig6}
 \epsfig{file=Loggamma.eps, height=200pt, width=160pt,angle=0}
                   }%
\caption{Estimations of $1/\alpha=1/2(a), 1/4(b), 1/4(c)$ (above)
and $\pk{X>3}=0.01298(a), 0.0035(b), 0.0037(c)$ (bottom), which is
indicated by the horizontal line. The line and the dotted line stand
for the estimators based on the original samples from $R, S$ and
$RS$, respectively. The true value of $\tau=-1(a), -2(b), 0(c)$.}
\end{center}
\end{figure}
\begin{figure}
\begin{center}
\subfigure[$R\sim |N(0,1)|,S\sim beta(1,3).$
]{%
            \label{fig7}
 \epsfig{file=Absolute(n).eps, height=200pt, width=160pt,angle=0}
                   }%
                    \subfigure[$R\sim W(2, 1), S\sim beta(2,3).$]{%
            \label{fig8}
 \epsfig{file=Weibull.eps, height=200pt, width=160pt,angle=0}
                   }%
                   \\
\subfigure[$R\sim \Gamma(1/2,1), S\sim beta(1,3).$]{%
            \label{fig9}
 \epsfig{file=Gamma(tp).eps, height=200pt, width=160pt,angle=0}
                   }%
  \subfigure[$R\sim PW(9/10,9/20), S\sim beta(1,3).$]{%
            \label{fig10}
 \epsfig{file=PWeibull.eps, height=200pt, width=160pt,angle=0}
                   }%
                   \caption{Sample paths of  \eee{estimators} of $\theta=1/2(a), 1/2(b),
1(c), 4/3(d)$ (above) and $\log(\widehat p_k/p)$ with
$\pk{X>3}=5.281\times 10^{-6}(a), 2.1186\times 10^{-7}(b),
\pk{X>10}=1.7189\times 10^{-8}(c), \lcx{3.5539\times 10^{-9}(d)}$
(bottom), which is indicated by the horizontal line. \lcx{The line
and the dotted line stand for the estimators based on the original
samples from $R, S$ and $RS$, respectively.} The true value of
$\tau_2=-\IF(a, c, d), -1(b)$ and $\rho'=-1(a), -\IF(b), -1(c),
-1/2(d)$.}
\end{center}
\end{figure}
}
\newpage

\subsection{Linear combinations of random contractions}
\EE{Motivated by the dependence structure of elliptical random vectors Hashorva et al. (2010) discussed the first-order tail asymptotics of the
aggregated risks of certain bivariate random vectors which we shall introduce next. Let therefore $(V_1,V_2)$ be
a bivariate scale mixture random vector  with stochastic representation}
\begin{equation}\label{4.0}
    (V_1,V_2)\stackrel{d}{=}R(I_1S, I_2\sqrt{1-S^2}),
\end{equation}
where $R\sim F,$ is almost surely positive, $S\sim G$ is a scaling
random variable taking values in $(0,1)$, while $I_1,I_2$ assume
values in $\{1,-1\}$.
%Furthermore, we suppose that $I_1, I_2, R$ and $S$ are mutually independent. \\
%The random vector $(V_1, V_2)$ is
%spherically distributed if $S^2\sim beta(1/2, 1/2)$ and
%$\pk{I_1=1}=\pk{I_2=1}=1/2$. Hence by Lemma 6.1 of Berman\,(1983)
%\begin{equation}\label{4.1}
%c_1V_1+c_2V_2\stackrel{d}{=}\sqrt{c_1^2+c_2^2}V_1,\quad \forall\
%c_1,c_2\in\R. \end{equation}
 %If $S\sim beta(a, b)$, then $(V_1, V_2)$ is a generalized Dirichlet random vector (cf. Hashorva et al.\,(2007)), and
 %\eqref{4.1} does not hold in general.
 %\\
 Hashorva et al.\,(2010) studied the
tail asymptotics of the aggregated risk
\begin{equation}\label{4.3}
    V(\eee{\lambda} )=\eee{\lambda}
V_1+\sqrt{1-\eee{\lambda} ^2}V_2=R( I_1\eee{\lambda}
S+I_2\sqrt{1-\eee{\lambda} ^2}\sqrt{1-S^2})=:RS(\eee{\lambda} )
\end{equation}
for \lcx{$\eee{\lambda} \in(0,1)$}.
\COM{It considered the following
assumption about the local \ck{behaviour} of $G$ at $\eee{\lambda} $
\begin{equation}\label{LLL}
\pk{|S-\eee{\lambda} |\le
x}=L_{\eee{\lambda} }(x)x^{\alpha_{\eee{\lambda} }}, \quad
    \alpha_{\eee{\lambda} }\in[0,\IF)
    \end{equation}
holds for all $x\in(0,\epsilon), \epsilon>0$, where
$L_{\eee{\lambda} }$ is \peng{a} positive slowly varying at $0$, and
$L_{\eee{\lambda} }(0+)=0$ if $\alpha_{\eee{\lambda} }=0$.}In what
follows, we derive the second-order tail asymptotics of
$V(\eee{\lambda} )$ in \eqref{4.3} if the following condition holds
for \eee{small} $x>0$
\begin{equation}\label{4.4}
 \pk{|S-\eee{\lambda} |\le
x}=c_{\eee{\lambda} } x^{\alpha_{\eee{\lambda} }}(1+L_{\eee{\lambda} }(x)x^{\tau_{\eee{\lambda} }}),
\quad
    \alpha_{\eee{\lambda} }, \tau_{\eee{\lambda} }\in(0,\IF)\quad\mbox{\lcx{and}}\quad \lambda\in[0,1],
    \end{equation}
    where $c_{\eee{\lambda} }$ is a positive constant and $|L_{\eee{\lambda} }|$ is
    slowly varying at $0$. Set
\begin{equation}\label{qlambda}
    q_\lambda= \pk{I_1=I_2=1}\I{\lambda\in(0,1)} +     \pk{I_2=1}\I{\lambda=0} +   \pk{I_1=1}\I{\lambda=1}.
          \end{equation}
          with $\I{\cdot}$ the indicator function.
\BL\label{L4.1}
    Let $I_1, I_2$ be two random variables taking values $-1,1$ with
    probability $q_\lambda\in(0,1]$ defined by \eqref{qlambda} and independent of the
    scaling random variable $S\sim G$. For
    given $\eee{\lambda} \in\lcx{[0,1]}$, if \lcx{further} df $G$ satisfies \eqref{4.4} for
    small $x>0$, then
    the  $S(\eee{\lambda})$ defined in \eqref{4.3} satisfies
    \COM{\begin{align}
       \nonumber \pk{S(\eee{\lambda} )>1-x}&=
        \lcx{q_\lambda}c_{\eee{\lambda} }\left[(2x(1-{\eee{\lambda} }^2))^{\alpha_{\eee{\lambda} }/2}\I{\lambda\in[0,1)} + x^{\alpha_\lambda}\I{\lambda=1}\right]\\
       \label{L4.5} &\quad\times \left[1-\left(\frac{\alpha_\lambda
       x}{4}\I{\lambda=0} + \frac{\alpha_{\eee{\lambda} }\eee{\lambda} }{\sqrt{2(1-{\eee{\lambda} }^2)}}x^{1/2}\I{\lambda\in(0,1)} \right)(1+o(1))\right.\\
       &\quad\qquad +\left.\left(L_{\eee{\lambda} }(\sqrt x)(2x(1-{\eee{\lambda} }^2))^{\tau_{\eee{\lambda}
       }/2}\I{\lambda \in[0,1)} + L_\lambda(x) x^{\tau_\lambda}\I{\lambda=1}
        \right)(1+o(1))\right],\quad x\downarrow0.
         \end{align}}
\begin{itemize}
  \item [a)] For $\lambda\in(0,1)$,
    \begin{equation*}%\label{L4.5}
        \pk{S(\eee{\lambda} )>1-x}\cH{=}
        q_\lambda c_{\eee{\lambda} }(2x(1-{\eee{\lambda} }^2))^{\alpha_{\eee{\lambda}
        }/2}\left[1+\mathcal{A}_\lambda(x)\right],
         \end{equation*} with
         $$\mathcal{A}_\lambda(x) = \left(L_{\eee{\lambda} }(\sqrt x)(2x(1-{\eee{\lambda} }^2))^{\tau_{\eee{\lambda} }/2} - \frac{\alpha_{\eee{\lambda} }\eee{\lambda} }{\sqrt{2(1-{\eee{\lambda}
         }^2)}}x^{1/2}\right)(1+o(1)).
        $$
  \item [b)] For $\lambda=0$,
  \begin{align*}
       %\label{L4.50}
       \pk{S(\eee{\lambda} )>1-x}&=
       q_\lambda c_{\eee{\lambda} }(2x)^{\alpha_{\eee{\lambda} }/2}\left[1 + \mathcal{A}_\lambda(x)\right],
       \quad \mathcal{A}_\lambda(x) = \left(L_{\eee{\lambda} }(\sqrt x)(2x)^{\tau_{\eee{\lambda}
       }/2} - \frac{\alpha_\lambda
       x}{4}\right)(1+o(1)).
         \end{align*}
  \item [c)] For $\lambda=1$,
  \begin{align*}
       %\label{L4.51}
       \pk{S(\eee{\lambda} )>1-x}&=
       q_\lambda c_{\eee{\lambda} }x^{\alpha_{\eee{\lambda} }}\left[1+\mathcal{A}_\lambda(x)\right],\quad \mathcal{A}_\lambda(x) = L_{\eee{\lambda} }(x)x^{\tau_{\eee{\lambda}
       }}.
         \end{align*}
\end{itemize}
\EL
\zz{In view of \nelem{L4.1} we have} $\pk{S(\eee{\lambda}
)>1-1/x}\in 2RV_{-\alpha,\tau}$ with $\alpha, \tau$ and auxiliary
function \cll{$A$ as follows}
\begin{equation}\label{A for
T4.1}
    \alpha=\left\{\begin{array}{ll}
                    \alpha_\lambda/2, & \lambda\in[0,1), \\
                    \alpha_\lambda, & \lambda=1;
                  \end{array}
                  \right.
    \quad
    \tau=\left\{\begin{array}{ll}
                    -\min(\tau_{\eee{\lambda} },1)/2, & \lambda\in (0,1), \\
                    -\min(\tau_\lambda, 2)/2, & \lambda=0,\\
                    -\tau_\lambda, & \lambda=1;
                  \end{array}
                  \right. \quad A(x) = \tau \mathcal{A}_\lambda(\lcx{1/x}).
\end{equation}
\ck{    Now, utilizing} \netheo{t2}, \netheo{t3} and \nelem{L4.1},
we \ck{derive} \peng{the following second-order tail asymptotics of
    $V(\eee{\lambda} )$}.
    \BT\label{T4.1}Let $V(\eee{\lambda} )$ be defined in \eqref{4.3} for $\lcx{\lambda\in[0,1]}$ and \cll{satisfying}\COM{\zz{satisfying further}} the conditions of \nelem{L4.1}.
\begin{enumerate}
  \item[a)] If $F\in D(Q_0)$ and its tail quantile
  function $U\in 2ERV_{0,\rho}, \cll{\rho\le0}$ with auxiliary functions $1/w(U)$ and $\tilde{A}$,
  then for $x\uparrow x_F$ (set $\lcx{\eta(x) =xw(x)}$)
\begin{align*}
    \nonumber \pk{V(\eee{\lambda} )>x} &= \bar F(x) \mathbb{P}\left(S(\lambda)>1-\frac1{\eta(x)}\right)\\
     &\quad\times\left[\Gamma(\alpha+1)+
     \left(\frac{\Gamma(\alpha-\tau+1)-\Gamma(\alpha+1)}{\tau}
   A(\eee{\eta(x)})+K(\alpha,\rho) \tilde{A}\left(\frac 1{\bar
    F(x)}\right)\right)(1+o(1))\right].
\end{align*}
 \item[b)] If $F\in D(Q_{-1/\alpha_1}), \alpha_1>0$ and $x_F=1$. Furthermore, we assume its tail quantile function $U$ satisfies $1-U\in 2RV_{-1/\alpha_1,
 \tau_1/\alpha_1}$ with auxiliary function $\tilde{A}$, then for
 $x\downarrow0$
\begin{align*}
  \nonumber \pk{V(\eee{\lambda} )>1-x}& = \bar F(1-x)\pk{S(\lambda)>1-x}\\
 \nonumber &\quad \times  \left[\alpha_1B\left(\alpha_1,
  \alpha+1\right)  +\left(\frac{\alpha\alpha_1^2}{\tau_1}\left[B\left(\alpha,\alpha_1+1
  \right)-B\left(\alpha,\alpha_1-\tau_1+1\right)\right]
   \tilde{A}\left(\frac1{\bar F(1-x)}\right) \right.\right.\\
&\qquad\left.\left.\cl{+}\frac{\alpha_1}{\tau}\left[B\left(\alpha_1,\alpha-\tau+1\right)
-B\left(\alpha_1,\alpha+1\right)\right]A\left(\frac1{x}\right)\right)(1+o(1))\right].\end{align*}
\end{enumerate}
Here $ \alpha, \tau$ and $A$ are \peng{those} defined in \eqref{A
for T4.1}, and $\pk{S(\lambda)> 1-x}$ is given by \nelem{L4.1}. \ET

\EE{{\remark \label{rem4.3}
\begin{itemize}
     \item [a)] If $S$ has Beta distribution with positive parameters $a$ and $b$,
%If $S\sim beta(a, b)$ with $a,b>0$,
then \eqref{4.4} holds for $\lambda= 0, 1$ and $\alpha_0 =  a,
\alpha_1 = b, \tau_0 = \tau_1 = 1,$
$$ c_0 =\frac1{aB(a,b)},\quad L_0(x) =
-\frac{(b-1)a}{a+1}(1+o(1)),\quad c_1 = \frac1{b B(a,b)},\quad
L_1(x) = -\frac{(a-1)b}{b+1}(1+o(1)).$$
% Thus the results of \netheo{T4.1} for $a=b=1/2$ agree with Remark \ref{rem4.3} $a)$ and our main findings
% in Section \ref{sec3}.
\item [b)]
If $G$
% has a continuous pdf $g$, then condition \eqref{LLL} holds
% with $\alpha_\lambda=1, L_\lambda(x) = ( 2+o(1)) g(\lambda)$ for all
% $\lambda\in (0,1)$. If \lcx{further} $G$
 has a continuous \lcx{3rd}
differentiable pdf $g$, then condition \eqref{4.4} holds with all
\lcx{$\lambda\in(0,1)$} and
$$\alpha_\lambda = 1,\quad c_\lambda = 2g(\lambda), \quad L_\lambda(x) = \frac{g\lcx{'''}(\lambda)}{6g'(\lambda)}\lcx{(1+o(1))},\quad \tau_\lambda = 2.$$
\item [c)] If $S$ has Beta distribution  with parameters $1/2, 1/2$
% If $S\sim beta(1/2,1/2)$
     and $I_1,I_2$ are independent with mean 0 being further
independent of $S$, then $(V_1,V_2)$ is spherically distributed.
\lcx{And} $V(\lambda) \stackrel{d}{=}I_1 RS \stackrel{ d}{=}I_2
R\sqrt{1-S^2}$ for all $\lambda\in[0,1]$. \lcx{Thus} the tail
asymptotics of $\lcx{V}(\lambda)$ can be directly obtained by
\lcx{\netheo{t2} and \netheo{t3}} in Section \ref{sec3}.
\end{itemize}}}

 \COM{\EE{Motivated by the dependence structure of elliptical random vectors Hashorva et al. (2010) discussed the first-order tail asymptotics of the
aggregated risks of certain bivariate random vectors which we shall
introduce next. Let therefore $(V_1,V_2)$ be a bivariate scale
mixture random vector  with stochastic representation}
\begin{equation}\label{4.0}
    (V_1,V_2)\stackrel{d}{=}R(I_1S, I_2\sqrt{1-S^2}),
\end{equation}
where $R\sim F,$ is almost surely positive, $S\sim G$ is a scaling
random variable taking values in $(0,1)$, while $I_1,I_2$ assume
values in $\{1,-1\}$.
%Furthermore, we suppose that $I_1, I_2, R$ and $S$ are mutually independent. \\
%The random vector $(V_1, V_2)$ is
%spherically distributed if $S^2\sim beta(1/2, 1/2)$ and
%$\pk{I_1=1}=\pk{I_2=1}=1/2$. Hence by Lemma 6.1 of Berman\,(1983)
%\begin{equation}\label{4.1}
%c_1V_1+c_2V_2\stackrel{d}{=}\sqrt{c_1^2+c_2^2}V_1,\quad \forall\
%c_1,c_2\in\R. \end{equation}
 %If $S\sim beta(a, b)$, then $(V_1, V_2)$ is a generalized Dirichlet random vector (cf. Hashorva et al.\,(2007)), and
 %\eqref{4.1} does not hold in general.
 %\\
 Hashorva et al.\,(2010) studied the
tail asymptotics of the aggregated risk
\begin{equation}\label{4.3}
    V(\eee{\lambda} )=\eee{\lambda}
V_1+\sqrt{1-\eee{\lambda} ^2}V_2=R( I_1\eee{\lambda}
S+I_2\sqrt{1-\eee{\lambda} ^2}\sqrt{1-S^2})=:RS(\eee{\lambda}
),\quad \eee{\lambda} \in(0,1)
\end{equation}
 if the following assumption
about the local \ck{behaviour} of $G$ at $\eee{\lambda} $
\begin{equation}\label{LLL}
\pk{|S-\eee{\lambda} |\le x}=L_{\eee{\lambda}
}(x)x^{\alpha_{\eee{\lambda} }}, \quad
    \alpha_{\eee{\lambda} }\in[0,\IF)
    \end{equation}
holds for all $x\in(0,\epsilon), \epsilon>0$, where
$L_{\eee{\lambda} }$ is \peng{a} positive slowly varying at $0$, and
$L_{\eee{\lambda} }(0+)=0$ if $\alpha_{\eee{\lambda} }=0$. In what
follows, we derive second-order tail asymptotics of $V(\eee{\lambda}
)$ in \eqref{4.3} if the following condition holds for \eee{small}
$x>0$
\begin{equation}\label{4.4}
 \pk{|S-\eee{\lambda} |\le
x}=c_{\eee{\lambda} } x^{\alpha_{\eee{\lambda} }}(1+L_{\eee{\lambda}
}(x)x^{\tau_{\eee{\lambda} }}), \quad
    \alpha_{\eee{\lambda} }, \tau_{\eee{\lambda} }\in(0,\IF),
    \end{equation}
    where $c_{\eee{\lambda} }$ is a positive constant and $|L_{\eee{\lambda} }|$ is
    slowly varying at $0$.
\BL\label{L4.1}
    Let $I_1, I_2$ be two random variables taking values $-1,1$ with
    probability $q=\pk{I_1=I_2=1}\in(0,1]$ and independent of the
    scaling random variable $S\sim G$. For
    given $\eee{\lambda} \in(0,1)$, if df $G$ satisfies \eqref{4.4}, then
    for the  $S(\eee{\lambda} )$ defined in \eqref{4.3},
    we have
    \begin{equation}\label{L4.5}
        \pk{S(\eee{\lambda} )>1-x}\cH{=}
        qc_{\eee{\lambda} }(2x(1-{\eee{\lambda} }^2))^{\alpha_{\eee{\lambda} }/2}\left[1-\left(\frac{\alpha_{\eee{\lambda} }\eee{\lambda} }{\sqrt{2(1-{\eee{\lambda} }^2)}}x^{1/2}-L_{\eee{\lambda} }(\sqrt x)(2x(1-{\eee{\lambda} }^2))^{\tau_{\eee{\lambda} }/2}
        \right)(1+o(1))\right],\quad x\downarrow0.
         \end{equation}
\EL \zz{In view of \nelem{L4.1} we have} $\pk{S(\eee{\lambda}
)>1-1/x}\in 2RV_{-\alpha,\tau}$ with auxiliary function \cll{$A$ as
follows}
   \begin{equation}\label{A for T4.1}
    \alpha=\alpha_{\eee{\lambda} }/2,\quad
    \tau=-\min(\tau_{\eee{\lambda} },1)/2,\quad A(x)=\left\{\begin{array}{ll}
                    \frac{2\alpha_{\eee{\lambda} }{\eee{\lambda} }}{\sqrt{2(1-{\eee{\lambda} }^2)}}x^{-1/2}, & \tau_{\eee{\lambda} }>1, \\
                    -\frac{\tau_{\eee{\lambda} }}{2}L_{\eee{\lambda} }(x^{-1/2})(2(1-\eee{\lambda} ^2))^{\tau_{\eee{\lambda} }/2}x^{-\tau_{\eee{\lambda} }/2}, & \tau_{\eee{\lambda} }<1,\\
                    \left[\frac{2\alpha_{\eee{\lambda} }{\eee{\lambda} }}{\sqrt{2(1-{\eee{\lambda} }^2)}}-\frac{\sqrt{2(1-\eee{\lambda} ^2)}}{2}L_{\eee{\lambda} }(x^{-1/2})\right]x^{-1/2}, &
                    \tau_{\eee{\lambda} }=1.
                  \end{array}
    \right.
   \end{equation}
\ck{    Now, utilizing} \netheo{t2}, \netheo{t3} and \nelem{L4.1},
we \ck{derive} \peng{the following second-order tail asymptotics of
    $V(\eee{\lambda} )$}\cll{.}
    \BT\label{T4.1}Let $V(\eee{\lambda} )$ be defined in \eqref{4.3} \zz{satisfying further} the conditions of \nelem{L4.1}.
\begin{enumerate}
  \item[(a)] If $F\in D(Q_0)$ and its tail quantile
  function $U\in 2ERV_{0,\rho}, \cll{\rho\le0}$ with auxiliary functions $1/w(U)$ and $\tilde{A}$,
  then for $x\uparrow x_F$
\begin{align}\label{Result 4.1}
    \nonumber \pk{V(\eee{\lambda} )>x}&\cH{=} q\bar
    F(x)c_{\eee{\lambda} }\left(\frac{2(1-{\eee{\lambda} }^2)}{\eee{\eta(x)}}\right)^{\alpha_{\eee{\lambda} }/2}\left[1-\left(\frac{\alpha_{\eee{\lambda} }{\eee{\lambda} }}{\sqrt{2(1-{\eee{\lambda} }^2)\eee{\eta(x)}}}-\frac{L_{\eee{\lambda} }((\eee{\eta(x)})^{-1/2})}
    {(2(1-{\eee{\lambda} }^2)\eee{\eta(x)})^{\tau_{\eee{\lambda} }/2}}\right)(1+o(1))\right]\\
    &\quad \times \left[\Gamma(\alpha+1)+
     \left(\frac{\Gamma(\alpha-\tau+1)-\Gamma(\alpha+1)}{\tau}
   A(\eee{\eta(x)})+K(\alpha,\rho) \tilde{A}\left(\frac 1{\bar
    F(x)}\right)\right)(1+o(1))\right].
\end{align}
 \item[(b)] If $F\in D(Q_{-1/\alpha_1}), \alpha_1>0$ and $x_F=1$. Furthermore, we assume its tail quantile function $U$ satisfies $1-U\in 2RV_{-1/\alpha_1,
 \tau_1/\alpha_1}$ with auxiliary function $\tilde{A}$, then for
 $x\downarrow0$
\begin{align}
  \nonumber \pk{V(\eee{\lambda} )>1-x}&\cH{=} q\bar F(1-x)c_{\eee{\lambda} }(2x(1-\eee{\lambda} ^2))^{\alpha_{\eee{\lambda} }/2}\left[1-\left(\frac{\alpha_{\eee{\lambda} }\eee{\lambda} }{\sqrt{2(1-\eee{\lambda} ^2)}}x^{1/2}-L_{\eee{\lambda} }(\sqrt
  x)(2x(1-\eee{\lambda} ^2))^{\tau_{\eee{\lambda} }/2}\right)(1+o(1))
        \right]\\
        &\quad \times \left[\alpha_1B\left(\alpha_1, \alpha+1\right)
  +\left(\frac{\alpha\alpha_1^2}{\tau_1}\left[B\left(\alpha,\alpha_1+1
  \right)-B\left(\alpha,\alpha_1-\tau_1+1\right)\right]\right.
   \tilde{A}\left(\frac1{\bar F(1-x)}\right)\right.\nonumber\\
\nonumber
&\quad\quad\left.\left.\cl{+}\frac{\alpha_1}{\tau}\left[B\left(\alpha_1,\alpha-\tau+1\right)
-B\left(\alpha_1,\alpha+1\right)\right]A\left(\frac1{x}\right)\right)(1+o(1))\right].\end{align}
\end{enumerate}
Here $\tau, \alpha$ and $A$ are \peng{those} defined in \eqref{A for
T4.1}. \ET

\EE{{\bf Remarks}: a) If $S$ has Beta distribution with parameters $1/2, 1/2$ and $I_1,I_2$ are independent with mean 0 being further independent of $S$, then $(V_1,V_2)$ is spherically distributed.\\
b) If $G$ has a continuous pdf $g$, then condition \eqref{LLL} holds with ... Condition\\
c) If $S$ had Beta distribution with positive parameters $a$ and
$b$, then \eqref{LLL} holds with ... and \eqref{4.4} holds with. ...
In particular, when $(V_1,V_2)$ is a spherical random vector, so
$a=b=1/2$, then $V(\lambda) \equaldis V_1= I_1 R S$, and the result
of \netheo{} agrees with our main findings for the product $X=
RS$.}
}

\section{Proofs}\label{sec6}
 \prooftheo{t1} It follows from \peng{Breiman's Lemma} that
$$\lim_{x\to\IF}\frac{\bar H(x)}{\bar F(x)}=\E{S^{\alpha_1}}.$$
We \peng{consider two} cases $\tau_1<0$ and $\tau_1=0$ separately.
For $\tau_1<0$, by Lemma 5.2 of Draisma et al.\,(1999), for every
$\epsilon>0$, there exists $x_0=x_0(\epsilon)>0$ such that for all
$x>x_0$ and all $s\in(0,1)$
$$\left|\frac{\bar F(x/s)/\bar F(x)-s^{\alpha_1}}{\tilde{A}(x)}-s^{\alpha_1}\frac{s^{-\tau_1}-1}{\tau_1}\right|\le\epsilon(C_1+C_2s^{\alpha_1}+C_3s^{\alpha_1-\tau_1-\epsilon}),$$
with some positive constants $C_1, C_2$ and $C_3$ independent of $x$
and $s$. Therefore, by \peng{the} dominated convergence theorem
\begin{align*}
   \ &\lim_{x\to\IF}\frac{1}{\tilde{A}(x)}\left( \frac{\bar H(x)}{\bar F(x)}-\E{S^{\alpha_1}}\right)=\int_0^1\lim_{x\to\IF}\frac{\bar F(x/s)/\bar
   F(x)-s^{\alpha_1}}{\tilde{A}(x)}\,\d
   G(s)=\mathbb{E}\left\{S^{\alpha_1}\frac{S^{-\tau_1}-1}{\tau_1}\right\}.
\end{align*}
For $\tau_1=0$, note that for all $\alpha_1>0,$ the function
$f(s)=s^{\alpha_1}\log s^{-1}$ is continuous in $(0,1]$ and
$\lim_{s\to0+}f(s)=0$ implies that $f(s)$ is bounded in $[0,1]$ and
$\E{f(S)}$ exists. Similarly to the proof of the case $\tau_1<0$
\begin{align*}
   \ &\lim_{x\to\IF}\frac{1}{\tilde{A}(x)}\left( \frac{\bar H(x)}{\bar F(x)}-\E{S^{\alpha_1}}\right)=\mathbb{E}\left\{S^{\alpha_1}\log
   S^{-1}\right\}
\end{align*}
holds for the case $\tau_1=0$, hence the proof is complete. \QED

\prooftheo{t2} %First, we rewrite the left-hand side of \eqref{result2} as follows.
\cl{ Denote $t=1/\bar F(x)$}, \peng{noting that}
\begin{align*}
\bar H(x)&=\int_{x}^{x_F} \bar
    G\left(\frac{x}{y}\right)\, \d F(y)=\int_{t}^{\IF} \bar
    G\left(\frac{U(t)}{U(s)}\right)\, \d(1-1/s)=t^{-1}\int_{0}^{1} \bar
    G\left(1-\frac{U(t/s)-U(t)}{U(t/s)}\right)\, \d s
\end{align*}
\peng{and rewrite the left-hand side of \eqref{result2} as}
   \begin{align}\label{proof_1}
   \nonumber\frac{\bar H(x)}{\bar F(x)\bar G\left(1-\frac{1}{\eta(x)}\right)}&=
    \int _0^1
    \frac{ \bar G\left(1-\frac{U(t/s)-U(t)}{U(t/s)}\right)}{\bar G\left(1-\frac{a(t)}{U(t)}\right)}\,\d
    s\notag\\
     \nonumber &= \int _0^1
    \left(\frac{U(t/s)-U(t)}{a(t)}\frac{U(t)}{U(t/s)}\right)^{\alpha_2}\frac{L\left(\frac{U(t)}{a(t)}\middle\slash \left(\frac{U(t/s)-U(t)}{a(t)}\frac{U(t)}{U(t/s)}\right)\right)}{L\left(\frac{U(t)}{a(t)}\right)}\,\d
    s \notag\\
    &=\int_0^1( \PFF )^{\alpha_2}\frac{L(\cw{\Xi_t(s)})}{L(\varphi_t)}\,\d
    s,
   \end{align}
 where
% \PFF
$$ \PFF=\cH{q_t(s)\phi_t(s)} ,  \quad \Xi_t(s)=\frac{\varphi_t}{\PFF}, \quad \varphi_t=\frac{U(t)}{a(t)}$$
\peng{and}
$$ \PF=\frac{U(t/s)-U(t)}{a(t)},\quad a=1/w(U),\quad \phi_t(s)=\frac{U(t)}{U(t/s)}.$$
\peng{Further we decompose \eqref{proof_1} as}
\begin{align}
   \nonumber \frac{\bar H(x)}{\bar F(x)\bar G\left(1-\frac{1}{\eta(x)}\right)}-\Gamma(\alpha_2+1)&=\int_0^1\left((\PF)^{\alpha_2}-\log^{\alpha_2}(1/s)\right)\,\d s-\int_0^1(\PF)^{\alpha_2}(1-(\phi_t(s))^{\alpha_2})\,\d s
\\ &\quad+\int_0^1( \PFF )^{\alpha_2}\left(\frac{L(\Xi_t(s))}{L(\varphi_t)}-1\right)\,\d
    s=:I-II+III.\label{proof_2}
\end{align}
\cll{Since} \peng{\eqref{proof_1} tends to $\Gamma(\alpha_2+1)$ by
Theorem 3.1 in Hashorva et al.(2010). The rest is to derive the
convergence rates of the three terms on the right-hand side of
\eqref{proof_2}.}
 By Lemma 5.2 in Draisma et al.\,(1999), for every $\epsilon>0$,
there exists $t_0=t_0(\epsilon)>0$ such that for all $t>t_0$ and all
$s\in(0,1)$
$$\left|\frac{\PF-\log(1/s)}{\tilde{A}(t)}-\psi(1/s)\right|\le\epsilon(C_1+C_3s^{-\rho-\epsilon}),$$
with some positive constants $C_1$ and $C_3$, independent of $t$ and
$s$. Therefore, \peng{by Taylor's expansion and the dominated
convergence theorem, we have}
\begin{equation}\label{I}
\lim_{t\to\IF}\frac{I}{\tilde{A}(t)}=\int_0^1\alpha_2\log^{\alpha_2-1}(1/s)\psi(1/s)\,\d
s=K(\alpha_2,\rho),
\end{equation}
 with \peng{$\psi$ and $K(\alpha_2,\rho)$ defined in
\eqref{def2} and \eqref{Coefficient KK}, respectively}. \\
For the second term $II$, recall that $U\in \Pi(a)$ implies that
$U\in RV_0$ and $ \varphi_t\to\IF$ as $t\to\IF$. By Corollary B.2.10
of de Haan and Ferreira\,(2006), for all $s\in (0,1)$ and
sufficiently large $t$
\begin{equation}\label{Inequality_phi}
0\le\PF\le cs^{-\epsilon},\quad
0\le\phi_t(s)=\left(1+\frac{\PF}{\varphi_t}\right)^{-1}\le 1
\end{equation}
 for some
$c>1$ and any $\epsilon>0$. Hence,
 $$ \frac{1-\phi_t(s)}{1/\varphi_t}\le \PF\le
cs^{-\epsilon}.$$ Therefore by Taylor's expansion and \peng{the}
dominated convergence theorem
\begin{align}\label{II}
\lim_{t\to\IF}\frac{II}{1/\varphi_t}&=\alpha_2\int_0^1\log^{\alpha_2+1}(1/s)\,\d
s=\alpha_2\Gamma(\alpha_2+2).
 \end{align}
 Finally, we shall show below \eqref{III} holds for the third term
 $III$
\begin{align}\label{III}\nonumber&\lim_{t\to\IF}\frac{III}{A(\varphi_t)}-\frac{\Gamma(\alpha_2-\tau_2+1)-\Gamma(\alpha_2+1)}{\tau_2}\\
&\quad =\lim_{t\to\IF} \int_0^1( \PFF )^{\alpha_2}
\left(\frac{L(\Xi_t(s))/L(\varphi_t)-1}{A(\varphi_t)}-\frac{( \PFF
)^{-\tau_2}-1}{\tau_2}\right) \,\d s=0.
    \end{align}

Recall that $L\in2RV_{0,\tau_2}$ with auxiliary function $A$, by
Lemma 5.2 in Draisma et al.\,(1999), for every $\epsilon>0$, there
exists $t_0=t_0(\epsilon)>0$ such that for all $\varphi_t>t_0$, the
\cH{integral of the right-hand side} of \eqref{III} is dominated by
\begin{align}\nonumber
& \int_{\{s: s\in(0,1), \Xi_t(s)>t_0\}}\epsilon( \PFF
)^{\alpha_2}(C_1+C_3( \PFF )^{-\tau_2}\exp(\epsilon|\log( \PFF
)|))\,\d
s\\
\nonumber&\quad+ \int_{\{s:
s\in(0,1),
\Xi_t(s)<t_0\}}( \PFF )^{\alpha_2}\left|\frac{L(\Xi_t(s))/
L(\varphi_t)-1}{A(\varphi_t)}\right|\,\d
s\\
&\quad+\int_{\{s:
s\in(0,1),
\Xi_t(s)<t_0\}}( \PFF )^{\alpha_2}\left|\frac{( \PFF )^{-\tau_2}-1}{\tau_2}\right|\,\d
s=:J_1+J_2+J_3.\label{J}
\end{align}
Recall that \eqref{Inequality_phi} implies that
$f_t(s)=( \PFF )^{\alpha}, s\in(0,1)$ is integrable for all
$\alpha>0$ and sufficiently large $t$. Thus, $J_1$ tends to 0 since
$\epsilon$ is arbitrarily small, whereas $J_3$ tends to 0 due to
$\varphi_t/t_0\to\IF$. \\
\peng{To deal with $J_2$}, we need two inequalities of $L$ and $A$
stated below in \eqref{L_inequality} and \eqref{A_inequality}.
Indeed, note that $L\in 2RV_{0,\tau_2}, \tau_2<0$ implies that $L$
is ultimately bounded away from 0 and
$$L(t)=t^{\alpha_2}\bar G(1-1/t)\le t^{\alpha_2},\quad L(t)>1/M$$
hold for some given $M>0$ and sufficiently large $t$. By \cll{Potter
bounds} (cf. Proposition B.1.9 in de Haan and Ferreira\,(2006)), for
any $\epsilon>0$, there exists $t_0=t_0(\epsilon)>0$ such that
$\min(\varphi_t,\Xi_t(s))>t_0$
$$\frac{L(\Xi_t(s))}{L(\varphi_t)}\le c\max(( \PFF )^\epsilon,
( \PFF )^{-\epsilon}),$$ otherwise for
$\varphi_t>t_0,\Xi_t(s)\le t_0$ such that
\begin{equation}\label{L_inequality}
\frac{L(\Xi_t(s))}{L(\varphi_t)}\le\frac{(\Xi_t(s))^{\alpha_2}}{1/M}\le
Mt_0^{\alpha_2}.
\end{equation}
For $A$, note that $|A|\in RV_{\tau_2}$ and it is ultimately
decreasing. By \cl{the Karamata Representation (cf. Resnick\,(1987),
\pzx{p.17})}, for any given $\delta>0$ and
 $t_0<\varphi_t< \PFF t_0$
\begin{align}\label{A_inequality}
|A(\varphi_t)|\ge|A( \PFF t_0)|\ge
K_2( \PFF )^{\tau_2-\delta}|A(t_0)|,
\end{align}
 with some $K_2\in(0,1)$ a
constant. Therefore, the integrand of $J_2$ is dominated by
$$\frac{Mt_0^{\alpha_2}+1}{K_2|A(t_0)|}( \PFF )^{\alpha_2-\tau_2+\delta}\le \frac{Mt_0^{\alpha_2}+1}{K_2|A(t_0)|}(cs^{-\epsilon})^{\alpha_2-\tau_2+\delta}.$$
\peng{So, by the} dominated convergence theorem, $J_2$ tends to 0 as
$t\to\IF$. Thus this together with the proved results for $J_1$ and
$J_3$ concludes the proof of \eqref{J}, and thus \eqref{III} holds.
\netheo{t2} follows from \eqref{I}, \eqref{II} and \eqref{III}.
 \QED

\proofkorr{Corr1} \ck{For} $\peng{a=1/w(U)}$ the first-order
auxiliary function of $U$, it follows from Theorem B.3.1 in de Haan
and Ferreira\,(2006) that $a\in 2RV_{0,\rho}, \cl{\rho\le0}$ with
auxiliary function $\widetilde A$. Thus for sufficiently large $x$
\begin{equation}\label{Corr_1}
    \frac{w\left(x+\frac
    z{w(x)}\right)}{w(x)}\cl{=}1-\frac{e^{\rho z}-1}{\rho}\widetilde A\left(\frac1{\bar
    F(x)}\right)(1+o(1))
\end{equation}
holds for all $z\in\R$ (here \cll{$(e^{\rho z}-1)/{\rho}$ is
interpreted as $z$ for $\rho=0$}). Note that $\bar G(1-1/x)\in
2RV_{-\alpha_2,\tau_2}$ and $|A|\in RV_{\tau_2}$ yield that
\begin{align}\label{Corr_2}
\nonumber &\quad\frac{\bar G\left(1-\frac
1{\eta(x+z/w(x))}\right)}{\bar
G(1-1/\eta(x)}\cl{=}\left(\frac{\eta(x+z/w(x))}{\eta(x)}\right)^{-\alpha_2}\left(1+\frac{\left(\frac{\eta(x+z/w(x))}{\eta(x)}\right)^{\tau_2}-1}{\tau_2}A(\eta(x))(1+o(1))\right)\\
 &
=\left(\frac{x+z/w(x)}{x}\frac{w(x+z/w(x))}{w(x)}\right)^{-\alpha_2}\left[1+o(1/{\eta(x)})+o(\widetilde
A(1/\bar F(x)))\right]
=1-\left[\frac{\alpha_2z}{\eta(x)}-\alpha_2\frac{e^{\rho
z}-1}{\rho}\widetilde A\left(\frac1{\bar
    F(x)}\right)\right](1+o(1)).
\end{align}
Recall that $U\in 2ERV_{0,\rho}$ with auxiliary function $\widetilde
A$,
\begin{equation}\label{Corr_3}
\frac{\bar F\left(x+\frac z{w(x)}\right)}{\bar
F(x)}=e^{-z}\left(1+\psi(e^{-z})\widetilde{A}\left(\frac1{\bar
F(x)}\right)\right).
\end{equation}
The claim \eqref{Corr} follows from \eqref{result2}, \eqref{Corr_1},
\eqref{Corr_2}, \eqref{Corr_3} and \cH{the fact that
\begin{align}\label{A and eta}
\lim_{x\to\IF}\eta(x)\widetilde{A}\left(\frac1{\bar
F(x)}\right)=\lim_{x\to\IF}\frac{\widetilde{A}(x)}{a(x)/U(x)}=0
\end{align}
 for
$\rho<0$ (cf. Lemma B.3.16 in de Haan and
Ferreira\,(2006))}.\\% (notation: $\peng{a=1/w(U)}$)).} \\
\cH{Using \eqref{A and eta} and} $h(h^\leftarrow(t))\sim t$ for
$h=1/\bar H$ for \eqref{Corr}, one can verify that $U_X\in
2ERV_{0,0}$ with auxiliary functions stated by \eqref{CORR}.\QED

\proofkorr{Corr2} First, note that $\bar F(x)=\exp(-V(x))$ and
$$U(t)=V^\leftarrow(\log t)=(\log t)^\theta\ell(\log t)$$
and thus
\begin{align*}
U(tx)&=V^\leftarrow(\log tx)=(\log t)^\theta\ell(\log
t)\left(1+\frac{\log x}{\log t}\right)^\theta\frac{\ell(\log
t(1+\log
x/\log t))}{\ell(\log t)}\\
&\cl{=} U(t)\left(1+\theta\frac{\log x}{\log
t}+\frac{\theta(\theta-1)}{2}\frac{\log^2x}{\log^2t}(1+o(1))\right)\left(1+b(\log
t)\frac{(1+\log x/ \log t)^{\rho'}-1}{\rho'}(1+o(1))\right).
\end{align*}
Therefore, $U\in2ERV_{0,0}$ with auxiliary functions $a$ and
$\widetilde A$ as
\begin{align*}
    a(t)=\frac{\theta+b(\log t)}{\log t} U(t),\quad \widetilde
    A(t)=\frac{\theta-1}{\log t}.
\end{align*}
This implies that
\begin{equation}\label{WTD0}
\eta(x)=\frac{x}{a(1/\bar F(x))}=\frac{V(x)}{\theta+b(V(x))},\quad
\widetilde
    A\left(\frac1{\bar F(x)}\right)=\frac{\theta-1}{V(x)}.
\end{equation}
    By \netheo{t2},
\begin{align}
\bar H(x)
%&\cl{=}& \bar
 %       F(x)\bar G\left(1-\frac1{\eta(x)}\right)\Gamma(\alpha_2+1)\left[1+
 %    \left(\frac{\frac{\Gamma(\alpha_2-\tau_2+1)}{\Gamma(\alpha_2+1)}-1}{\tau_2}
 %  A(\eta(x))-
 %  \frac{\alpha_2(\alpha_2+1)}{\eta(x)}+\frac{\alpha_2(\alpha_2+1)}{2}\tilde{A}\left(\frac 1{\bar
 %   F(x)}\right)\right)(1+o(1)) \right]\\
 %
   \nonumber
    &= \bar F(x)\bar G\left(1-\frac1{V(x)}\right)\left(\frac{\eta(x)}{V(x)}\right)^{-\alpha_2}\left[1+\frac{
    \left(\frac{\eta(x)}{V(x)}\right)^{\tau_2}-1}{\tau_2}A(V(x))(1+o(1))\right]
    \Gamma(\alpha_2+1)\\
    \nonumber &\quad
\cl{\times}\left[1+\left(\frac{\displaystyle\frac{\Gamma(\alpha_2-\tau_2+1)}{\Gamma(\alpha_2+1)}-1}{\tau_2}\left(\frac{\eta(x)}{V(x)}\right)^{\tau_2}
A(V(x))-\frac{(\theta+b(V(x)))\alpha_2(\alpha_2+1)}{V(x)}
+\frac{(\theta-1)\alpha_2(\alpha_2+1)}{2V(x)}\right)(1+o(1))\right]\\
    \nonumber&=\exp(-V(x))\bar G\left(1-\frac1{V(x)}\right)\Gamma(\alpha_2+1)\theta^{\alpha_2}\\
  \label{WTD1} &\quad\times
\left[1+\left(\frac{\alpha_2}{\theta}b(V(x))+\frac{\displaystyle\frac{\Gamma(\alpha_2-\tau_2+1)}{\theta^{\tau_2}\Gamma(\alpha_2+1)}-1}{\tau_2}
A(V(x))-\frac{(\theta+1)\alpha_2(\alpha_2+1)}{2V(x)}\right)(1+o(1))\right]\\
     \label{WTD2}   &=:\exp(-V(x))(V(x))^{-\alpha_2}
        L^*(V(x)),
\end{align}
where \lcx{\eqref{WTD1}} is due to \eqref{WTD0} and $\bar
G(1-1/x)\in 2RV_{-\alpha_2,\tau_2}$ with auxiliary function $A$.
Clearly, $L^*$ is a slowly varying function. Therefore, let the
\lcx{right-hand} side of \eqref{WTD2} equal to $1/s$, and solve the
equation of $x$, then $V(x)\sim \log s$ and
\begin{align*}
    U_X(s)=V^\leftarrow\left(\log
    \frac{sL^*(V(x))}{(V(x))^{\alpha_2}}\right)&\lcx{=}\left(\log s-\alpha_2\log V(x)\left(1-\frac{\log L^*(V(x))}{\alpha_2\log V(x)}\right)\right)^\theta\ell\left(\log s-\alpha_2\log V(x)\left(1-\frac{\log L^*(V(x))}{\alpha_2\log V(x)}\right)\right)\\
    &=\left(\log s - \alpha_2\log\log s (1+o(1))\right)^\theta\ell(\log s)(1+o(\log\log s/\log s)).
\end{align*}
The last step is due to $\ell\in 2RV_{0,\rho'}$ and the property of
slowly varying function: $ \log L^*(V(x))/\log V(x)\to0$ (see
Bingham et al.\,(1987)). Hence
$$\bar H(x)=\exp(-V^*(x)),\quad (V^*)^\leftarrow(x)=x^\theta\left(1-\alpha_2\frac{\log x}{x}\right)^\theta\ell^*(x).$$
Thus the claim in \nekorr{Corr2} follows from $\ell^*\in
2RV_{0,\rho'^*}$ with $\rho'^*=\max(\rho', -1)$ and auxiliary
function $$b^*(x)=b(x)+\frac{\theta\alpha_2\log x}x.$$\QED

 \prooftheo{t3} First, by arguments similar to the case $F\in D(Q_{0})$ (cf. \eqref{proof_1}\cll{)}, we have
\begin{align*}
   \nonumber&\quad \frac{\bar H(x)}{\bar F(x)\bar
G(x)}
%=\int_0^1 \frac{\bar
%    G\left(1-\frac{U(t/s)-U(t)}{U(t/s)}\right)}{\bar
%    G(1-(1-U(t))}\, \d s\\
%   \nonumber &=\int_0^1\left(\frac{U(t/s)-U(t)}{1-U(t))}\frac{1}{U(t/s)}\right)^{\alpha_2}\frac{L\left(\frac{1}{1-U(t)}\middle\slash %\left(\frac{U(t/s)-U(t)}{1-U(t)}\frac{1}{U(t/s)}\right)\right)}{L\left(\frac{1}{1-U(t)}\right)}\,\d
%    s\\
 %   &
 =\int_0^1( \PFF )^{\alpha_2}\frac{L(\frac{\varphi_t}{ \PFF })}{L(\varphi_t)}\,\d
    s,
\end{align*}
with $t=1/\bar F(x), x=U(t)$ and
\begin{equation*}\label{E1}
 \PFF= \PF \phi_t(s),\quad \varphi_t=\frac{1}{1-U(t)}
 \quad \mbox{\lcx{with}} \quad\phi_t(s)=\frac{1}{U(t/s)},
\quad \PF=\frac{U(t/s)-U(t)}{1-U(t)}.
\end{equation*}
 \peng{Next}
% \xx{Consequently, we have }
\begin{align}
   \nonumber \frac{\bar H(x)}{\bar F(x)\bar G(x)}-\alpha_1B(\alpha_1,\alpha_2+1)&=\int_0^1\big(\PF\big)^{\alpha_2}-(1-s^{1/\alpha_1})^{\alpha_2}\,\d s\\
   \nonumber&\quad+\int_0^1(\PF)^{\alpha_2}((\phi_t(s))^{\alpha_2}-1)\,\d s
 +\int_0^1( \PFF )^{\alpha_2}\left(\frac{L(\frac{\varphi_t}{ \PFF })}{L(\varphi_t)}-1\right)\,\d
    s\\
    &=:I+II+III.\label{Proof_2}
\end{align}
\peng{\ck{It remains thus to} derive the convergence rate of each term
in \eqref{Proof_2}.} By Lemma 5.2 in Draisma et al.\,(1999), for
every $\epsilon>0$, there exists $t_0=t_0(\epsilon)>0$ such that for
all $t>t_0$ and all $s\in(0,1)$
$$\left|\frac{\PF-(1-s^{1/\alpha_1})}{\tilde{A}(t)}+s^{1/\alpha_1}\frac{s^{-\tau_1/\alpha_1}-1}{\tau_1/\alpha_1}\right|
\le\epsilon(C_1+C_2s^{1/\alpha_1}+C_3s^{(1-\tau_1)/\alpha_1\aH{-}\epsilon}),$$
with some positive constants $C_1, C_2$ and $C_3$, independent of
$t$ and $s$. Therefore, by Taylor's expansion and \peng{the}
dominated convergence theorem
\begin{equation}\label{3I}
    \lim_{t\to\IF}\frac{I}{\tilde{A}(t)}=-\alpha_2\int_0^1(1-s^{1/\alpha_1})^{\alpha_2-1}s^{1/\alpha_1}\frac{s^{-\tau_1/\alpha_1}-1}{\tau_1/\alpha_1}\,\d
    s=-\frac{\alpha_2\alpha_1^2}{\tau_1}(B(\alpha_2, \alpha_1-\tau_1+1)-B(\alpha_2,
    \alpha_1+1)).
\end{equation}
\cH{Here, \eqref{3I} for $\tau_1=0$ is understood \xx{as}
$$-\alpha_2\int_0^1(1-s^{1/\alpha_1})^{\alpha_2-1}s^{1/\alpha_1}\lim_{\tau_1\to0}\frac{s^{-\tau_1/\alpha_1}-1}{\tau_1/\alpha_1}\,\d
    s=\lim_{\tau_1\to0}-\frac{\alpha_2\alpha_1^2}{\tau_1}\big(B(\alpha_2, \alpha_1-\tau_1+1)-B(\alpha_2,
    \alpha_1+1)\big)$$
(cf. Corollary 4.4 in Mao and Hu\,(\lcx{2012(a)})).} For $II$, note
that $\PF\in(0,1), \varphi_t\to\IF$ and thus for all $s\in(0,1)$
$$0\le \frac{\phi_t(s)-1}{1/\varphi_t}=\frac{ \left(1-(1-\PF)/{\varphi_t}\right)^{-1}-1}{1/\varphi_t}=\frac{1-\PF}{1-(1-\PF)/{\varphi_t}}\le\frac1{1-1/\varphi_t}\to1$$
as $t\to\IF$. So, by Taylor's expansion and \peng{the} dominated
convergence theorem
\begin{align}\label{3II}
 \nonumber\lim_{t\to\IF}\frac{II}{1/\varphi_t}&\cll{=}\int_0^1\lim_{t\to\IF}(\PF)^{\alpha_2}\frac{\left(1+(\phi_t(s)-1)\right)^{\alpha_2}-1}{1/\varphi_t}\,\d s\\
 &=\alpha_2\int_0^1(1-s^{1/\alpha_1})^{\alpha_2}s^{1/\alpha_1}\,\d
    s=\alpha_1\alpha_2B(\alpha_1+1,\alpha_2+1).
\end{align}
Now we consider the third term $III$. By Lemma 5.2 in Draisma et
al.\,(1999), for every $\epsilon>0$, there exists
$t_0=t_0(\epsilon)>0$ such that for all $\varphi_t>t_0$ and all
$s\in(0,1)$
\begin{align*}
&\quad\left|( \PFF )^{\alpha_2}
\left(\frac{L(\frac{\varphi_t}{ \PFF })/L(\varphi_t)-1}{A(\varphi_t)}-\frac{( \PFF )^{-\tau_2}-1}{\tau_2}\right)\right|\\
    &\quad\quad\le
    \epsilon(C_1+C_2( \PFF )^{\alpha_2}+C_3( \PFF )^{\alpha_2-\tau_2-\epsilon}\cH{)}\le\epsilon(C_1+C_2+C_3).
\end{align*}
The last step is due to $ \PFF \le1$ for all $s\in(0,1)$ and $t>0$.
\peng{Hence, by the} dominated convergence theorem
\begin{align}\label{3III}
   \nonumber\lim_{t\to\IF} \frac{III}{A(t)}&\cll{=}\int_0^1\lim_{t\to\IF}( \PFF )^{\alpha_2}\frac{( \PFF )^{-\tau_2}-1}{\tau_2}\,\d s\\
   &=\int_0^1 (1-s^{1/\alpha_1})^{\alpha_2}\frac{(1-s^{1/\alpha_1})^{-\tau_2}-1}{\tau_2}\,\d
   s=
   \frac{\alpha_1}{\tau_2}\big(B\left(\alpha_1,\alpha_2-\tau_2+1\right)
-B\left(\alpha_1,\alpha_2+1\right)\big).
\end{align}
 The
claim follows from \eqref{3I}, \eqref{3II} and \eqref{3III}.
 \QED
 \\
\prooflem{L4.1} \COM{It follows that  $S(\eee{\lambda} )\le1$ and it
is bounded away from unity unless $I_1=I_2=1$ , and when the event
$\{I_1=I_2=1\}$ occurs, $S(\eee{\lambda} )\uparrow1$ if and only if
$|S-\eee{\lambda} |\downarrow0$. For small $x>0$, \cll{the event}
$$\{S(\eee{\lambda} )>1-x\}=\{(S-\eee{\lambda} )^2+2\eee{\lambda}  xS<2x-x^2\}$$
\cll{is} equivalent to $$(S-\eee{\lambda}
)^2<2x\big((1-\eee{\lambda} ^2)-\eee{\lambda}
\sqrt{2x(1-\eee{\lambda} ^2)}(1+o_p(1))\big).$$  Consequently, the
claim follows from \eqref{4.4}. \QED}
 We only give the proofs of the
case $\lambda\in(0,1)$. \lcx{The other cases} are left to the
readers and one can verify it by the similar arguments. Clearly, for
$\lambda\in(0,1), S(\eee{\lambda} )\le1$ and it is bounded away from
unity unless $I_1=I_2=1$, and when the event $\{I_1=I_2=1\}$ occurs,
$S(\eee{\lambda} )\uparrow1$ if and only if $|S-\eee{\lambda}
|\downarrow0$. For small $x>0$, \cll{the event}
$$\{S(\eee{\lambda} )>1-x\}=\{(S-\eee{\lambda} )^2+2\eee{\lambda}  xS<2x-x^2\}$$
\cll{is} equivalent to $$(S-\eee{\lambda}
)^2<2x\big((1-\eee{\lambda} ^2)-\eee{\lambda}
\sqrt{2x(1-\eee{\lambda} ^2)}(1+o_p(1))\big).$$
 Consequently, the
claim follows from \eqref{4.4}. \QED

\COM{ For $\lambda=0, S(0)=I_2\sqrt{1-S^2}$ and $S(0)$ is bounded
away 1 unless $I_2=1$. When $\I_2=1$, $\{S(0)>1-x\} = \{S^2<2x-x^2\}
= \{|S|<\sqrt{2x}(1-x/4(1+o_p(1)))\}$,
\\
For $\lambda=1, S(1)=I_1S$ and when $I_1=1, \{S(1)>1-x\} =
\{|S-1|<x\}$.}
\newpage

\section{Appendix}\label{appendix}
\COM{\peng{\bf A. \prooftheo{t_RP}}
 \prooftheo{t_RP}By Tang and
Vernic\,(2010), the ruin probability
$p(x)\sim\pk{RS>x}=\pk{RS/x_S>x_1}, x_1=x/x_S.$ Note that
\begin{equation}\label{3.01}\pk{S/x_S>1-1/x}=\pk{\Upsilon/(\pi x_S)+1>x} \end{equation} and
 $\pk{\Upsilon>x}\in 2RV_{-\alpha,\tau}$ with auxiliary
function $A$ satisfying
\begin{equation}\label{limit a}\lim_{t\to \IF}tA(t)=c\in[-\IF,\IF],
\end{equation} implying $\pk{\Upsilon/(\pi x_S)+1>x}\in
2RV_{-\alpha,\tau^*}$ with $\tau^*=\max(\tau,-1)$ and auxiliary
function
$$A^*(x)=\left\{
                  \begin{array}{ll}
                    -\alpha/x, & c=0, \\
                    (\pi x_S)^\tau A(x), & c=\pm\IF, \\
                    \frac{-\alpha+c/(\pi x_S)}x, & \mbox{otherwise}. \\
                  \end{array}
\right.$$Hence, by \netheo{t2}, we have
\begin{align*}
    \pk{RS>x}&=\pk{R{S}/{x_S}>x_1}\cl{=} \bar
    F(x_1)\mathbb{P}\left(\frac S{x_S}>1-\frac1{\eee{\eta(x_1)} }\right)\left[\Gamma(\alpha+1)+\left(\frac{\Gamma(\alpha-\tau^*+1)-\Gamma(\alpha+1)}{\tau^*}
   A^*(\eee{\eta(x_1)} )\right.\right.\\
    &\quad-\left.\left.
     \alpha\Gamma(\alpha+2)
   \frac 1{\eee{\eta(x_1)} }+K(\alpha,\rho)\tilde{A}\left(\frac 1{\bar
    F(x_1)}\right) \right)(1+o(1))\right]\\
    &\cl{=} \bar F(x_1)\pk{\Upsilon>\eee{\eta(x_1)} }\frac1{(\pi x_S)^\alpha}\left[1+\frac{(\pi x_S)^\tau-1}{\tau}A(\eee{\eta(x_1)} )(1+o(1))\right]\\
    &\quad\cH{\times}\left[\Gamma(\alpha+1)+\left(\frac{\Gamma(\alpha-\tau^*+1)-\Gamma(\alpha+1)}{\tau^*}
   A^*(\eee{\eta(x_1)} )-
        \frac{\alpha\Gamma(\alpha+2)}{\eee{\eta(x_1)} }+K(\alpha,\rho)\tilde{A}\left(\frac 1{\bar
    F(x_1)}\right) \right)(1+o(1))\right],
\end{align*}
where the last step is due to \eqref{3.01} and $\pk{\Upsilon>x}\in
2RV_{-\alpha,\tau}$ with auxiliary function $A$ and $|A|\in RV_\tau,
\tau<0$. This concludes the proof of \netheo{t_RP}. \QED}

\cll{This appendix includes two tables. Table \ref{T2} \ck{contains}
Weibull tail distributions satisfying the second-order regular
varying conditions and Table \ref{T1} shows several distributions in
maximum domain attraction of the \FRE distribution, the Gumbel
distribution and the Weibull distribution in the second-order
framework.}

\begin{minipage}{\linewidth}
\centering \captionof{table}{\peng{{\bf Weibull tail
distributions}}} \label{T2}
\begin{tabular}{l|l|c|c|c}
  \toprule[1.25pt]
  Weibull tail distributions & Tail $\bar F(x)$ \mbox{or\ pdf}\ $f(x)$ & $\theta$ & $\rho$ & $b(x)$ \\
  \midrule[0.8pt]
  Gamma ($\Gamma(\alpha, \lambda)$) & $f(x)=\frac{\lambda^\alpha}{\Gamma(\alpha)}x^{\alpha-1}e^{-\lambda x},
 \lambda, \alpha>0,\alpha\neq1 $ & 1 & $-1$ & $(1-\alpha){\log x}/{x}$
 \\
  \cll{Absolute} Normal ($|N(0,1)|$) & $f(x)=\frac{2}{\sqrt{2\pi}}e^{-x^2/2}$ & 1/2 &
  $-1$ & ${\log x}/({4x})$\\
     Weibull ($W(\beta, c)$)& $\bar F(x)=\exp(-c x^{\beta}), c,
  \beta>0$ & $1/\beta$ & $-\IF$ & 0\\
  Pertured Weibull ($PW(\beta,\alpha)$) & $\bar F(x)=e^{-x^\beta(C+Dx^{-\alpha})},
  \alpha,\beta, C>0, D\in \R$ & $1/\beta$ & $-\alpha/\beta$ & $\frac{\alpha D}{\beta^2}C^{\alpha/\beta-1}
  x^{-\alpha/\beta}$\\
  Modified Weibull ($MW(\beta)$) & $Y\log Y\sim F, Y\sim W(\beta)$ & $1/\beta$ & 0 & $1/\log
  x$\\
 Benktander II ($\mathcal{B}II(\beta, \lambda)$) & $\bar F(x)=x^{-(1-\beta)}\exp(-\frac\lambda\beta (x^\beta-1)), \lambda>0,
 0<\beta<1$ & $1/\beta$ & $-1$ & ${(1-\beta)\log x}/(\beta^2x)$ \\
  Extended Weibull ($\mathcal{E}W(\beta, \alpha)$) & $\bar F(x)=r(x)\exp(-x^\beta), \beta\in(0,1), r\in
 RV_{-\alpha}, \alpha\in \R$ & $1/\beta$ & $-1$ & $\alpha\log x/(\beta^2x)$ \\
 Logistic & $\bar F(x)=\frac2{1+e^x}$ & 1 & $-1$ &
 $-(\log2)/{x}$\\
 Gumbel ($G(\mu)$) & $\bar F(x)=1-\exp(-\exp(\mu-x)),
 \mu\neq0$ & 1 & $-1$ & $-\mu /x$\\
 \bottomrule[1.25pt]
\end{tabular}
\par
\bigskip
Weibull tail distribution: $\bar F(x)=\exp(-V(x)),
V^\leftarrow(x)=x^\theta \ell(x)$ and $\ell\in 2RV_{0,\rho}$ with
auxiliary function $b$.
\end{minipage}

%\begin{sidewaystable}[htbp]
\begin{minipage}{\linewidth}
\centering \captionof{table}{\peng{{\bf Risks satisfying the
second-order regular variation conditions}}} \label{T1}
\begin{tabular}{l|l|c|c|c}
  \toprule[1.25pt]
$\FRE$ attraction & Tail $\bar F(x)$ \mbox{or\ pdf}\ $f(x)$ & $\alpha$ & $\tau$ & $A(x)$ \\
 \midrule[0.8pt]
 Pareto& $\bar F(x)=\left(\frac{\theta}{\theta+x}\right)^{\alpha},\quad \theta,\alpha>0$& $\alpha$ &
 $-1$
 & $\alpha\theta/{x}$\\
 \FRE & $\bar F(x)=1-\exp(-x^{-\alpha})$ & $\alpha$ & $-\alpha$ &
$\alpha x^{-\alpha}/{2}$\\
 Burr& $\bar F(x)=(1+x^b)^{-a}$ & $ab$ & $-b$ & $ab x^{-b}$\\
 Hall-Weiss & $\bar F(x)=\frac12x^{-\alpha}(1+x^\tau),\quad
 \alpha>0,\tau<0$ & $\alpha$ & $\tau$ & $\tau x^\tau$\\
  $F(m,n)$ & $f(x)=\frac1{B(m/2,n/2)}\left(\frac mn\right)^{m/2}x^{m/2-1}\left(1+\frac{mx}n\right)^{-(m+n)/2}$&
 $n/2$& $-1$ &\cL{$\frac{(m+n)n^2}{2m(n+2)x}$}\\
 Log-gamma & $f(x)=\frac{\alpha^\beta}{\Gamma(\beta)}(\log
 x)^{\beta-1}x^{-\alpha-1}, \quad \alpha, \beta>0$ & $\alpha$ & 0 &
$ {(\beta-1)}/{\log x}$\\
Inv\cll{-}gamma &
$f(x)=\frac{\beta^\alpha}{\Gamma(\alpha)}x^{-\alpha-1}e^{-\beta/x},
\quad \alpha, \beta>0$ & $\alpha$ & $-1$ &
\cL{$ \frac{\alpha\beta}{(\alpha+1)x}$}\\
 Absolute \cll{$t$}&
$f(x)=\frac{2\Gamma(v/2)}{\sqrt{v\pi}\Gamma((v+1)/2)}(1+x^2\pzx{/v})^{-(v+1)/2},\quad
v\in \mathbb{N}$
& $v$ & $-2$ & $ \lcx{\frac{v^2(v+1)}{(v+2)x^2}}$\\
 \midrule[0.8pt]
  Weibull attraction & Tail $\bar F(x_F-1/x)$ \mbox{or\ pdf}\ $f(x)$ & $\alpha$ & $\tau$ & $A(x)$ \\
 \midrule[0.8pt]
 Beta & $f(x)=\frac1{B(a,b)}x^{a-1}(1-x)^{b-1},\quad
a, b>0$& $b$ & $-1\ (a\neq1)$
 & $\frac{b(a-1)}{(b+1)x}$\\
 Reverse-Burr& $\bar F(x_F-1/x)=(1+x^b)^{-a}$ & $ab$ & $-b$ & $ab x^{-b}$\\
 Extreme value Weibull & $\bar F(x_F-1/x)=1-\exp(-x^{-\alpha})$ & $\alpha$ &
 $-\alpha$ & $\alpha x^{-\alpha}/2$ \\
\midrule[0.8pt]
  Gumbel attraction & Tail $\bar F(x)$ \mbox{or\ pdf}\ $f(x)$ & $\rho$ & $a(x)$ & $A(x)$ \\
  \midrule[0.8pt]
  Gamma & $f(x)=\frac{\lambda^\alpha}{\Gamma(\alpha)}x^{\alpha-1}e^{-\lambda x},
\quad \cl{\lambda, \alpha>0}$ & 0 & $\left(1+\frac{\alpha-1}{\log
x}\right)\big\slash \lambda$ & $(1-\alpha)/\log^2x$ \\
  Absolute Normal & $f(x)=\frac{2}{\sqrt{2\pi}}e^{-x^2/2}$ & 0 &
  $\frac{U_1(2x)}{2\log(2x)}$ & $-1/(2\log x)$\\
  Log-normal & $f(x)=\frac{1}{\sqrt{2\pi}x}\exp(-\frac{\log^2
  x}{2})$ & 0 & $\frac{\exp(U_1(x))}{\sqrt{2\log x}}$ & $1/{\sqrt{2\log x}}$ \\
   Logistic & $\bar F(x)=\frac2{1+e^x}$ & $-1$ & $1$ & $1/{(2x)}$\\
  Truncated Gumbel & $\bar
  F(x)=\frac{1-\exp(-e^{-x})}{1-e^{-1}}$ & $-1$ &1 &$({1-e^{-1}})/({2x})$ \\
  Exponential with finite $x_F$ & $\bar F(x)=\exp(-\frac{c}{x_F-x}+\frac{c}{x_F}),\quad c>0,
  x_F>0$ & 0 & $\frac{c}{(\log x+c/x_F)^2}$ & $-2/{\log x}$ \\
  Weibull & $\bar F(x)=\exp(-c x^{\beta}),\quad c>0,
  \beta\in(0,1)$ & 0 & $\frac{(\log x)^{1/\beta-1}}{\beta
  c^{1/\beta}}$ & ${(1/\beta-1)}/{\log x}$\\
  Benktander I & $\bar F(x)=\left(1+\frac{2\beta}\alpha\log x\right)\exp(-\beta\log^2x-(\alpha+1)\log
  x)$ & 0 & $\frac{U_2(x)}{2\sqrt{\beta\log x}}$ & ${1}/({2\sqrt{\beta\log
  x}})$\\
 Benktander II & $\bar F(x)=x^{-(1-\beta)}\exp(-\frac\alpha\beta (x^\beta-1)), \alpha>0,
  0<\beta<1$ & 0 & $a^*(x)$& $(1/\beta-1)/{\log x}$ \\
  \midrule[0.8pt]
 \multicolumn{5}{c}{$a^*(x)=\frac{1-(1-\beta)/(\beta(\alpha/\beta+\log x))}{\beta(\alpha/\beta+\log x)}U(x),
 \quad\quad\quad U(x)= \left(\frac\beta\alpha((\alpha/\beta+\log x)-(1-\beta)\log U(x))\right)^{1/\beta}$}\\
 \multicolumn{5}{c}{$U_1(x)=\sqrt{2\log x}-\frac{\log(4\pi\log x)}{2\sqrt{2\log x}},\quad \qquad U_2(x)=\exp\left(-\frac{\alpha+1}{2\beta}
 +\sqrt{\frac{\log x}{\beta}}+\frac{\log\log x+\log(4\beta/\alpha^2)+(\alpha+1)^2/(2\beta)}{4\sqrt{\beta\log x}}\right)$}\\
 \bottomrule[1.25pt]
\end{tabular}
\par
\bigskip
For \FRE attraction, $\bar F\in 2RV_{-\alpha,\tau}$ with auxiliary
function $A$. For Weibull attraction, $\bar F(x_F-1/x)\in
2RV_{-\alpha,\tau}$ with auxiliary function $A$ \cll{and a finite
upper endpoint $x_F$}. For Gumbel attraction, the tail quantile
function $U\in 2ERV_{0,\rho}$ with the first-order auxiliary
function $a$ and the second-order auxiliary function $A$.
\end{minipage}

\newpage

%\pzx{ \noindent {\bf Acknowledgements}~~The third author was
%supported by the National Natural Science Foundation of China grant
%no.11171275 and the Natural Science Foundation Project of CQ no.
%cstc2012jjA00029.}

 \bibliographystyle{plain}

\end{document}